\documentclass{amsart}

\usepackage[mathcal]{euscript}
\usepackage{mathrsfs}
\usepackage{amsfonts,amssymb,amsbsy,amsmath}
\usepackage{graphicx,color}
\usepackage{epsfig}
\usepackage[all,cmtip]{xy}

\makeatletter
\let\uppercasenonmath\@gobble

\makeatother

\newcommand{\N}{\mathbb{N}}

\newcommand{\RP}{{\mathbb{R}}P}

\newcommand{\CP}{{\mathbb{C}}P}
\newcommand{\R}{\mathbb{R}}

\newcommand{\C}{\mathbb{C}}
\newcommand{\Z}{\mathbb{Z}}
\newcommand{\M}{\mathit{Map}}

\newcommand{\ba}{\begin{array}}
\newcommand{\ea}{\end{array}}

\newcommand{\Sa}{T^2}
\newcommand{\cj}{\mathit{conj}}

\theoremstyle{plain}
\newtheorem{t.}{Theorem}[section]
\newtheorem{l.}[t.]{Lemma}
\newtheorem{p.}[t.]{Proposition}
\newtheorem{c.}[t.]{Corollary}

\theoremstyle{definition}
\newtheorem{d.}[t.]{Definition}

\newtheorem{r.}[t.]{Remark}

\newcommand{\ovbir}{\includegraphics [width=0.28cm]{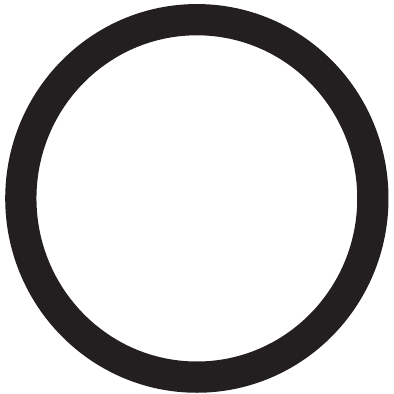}}
\newcommand{\oviki}{\includegraphics [width=0.28cm]{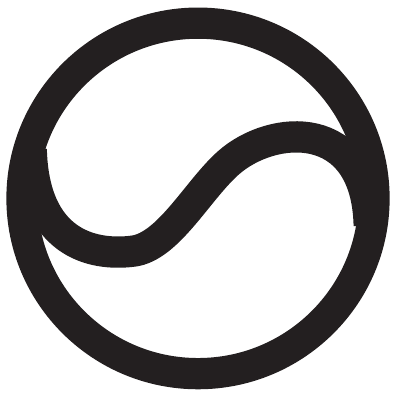}}
\newcommand{\ovuc}{\includegraphics [width=0.28cm]{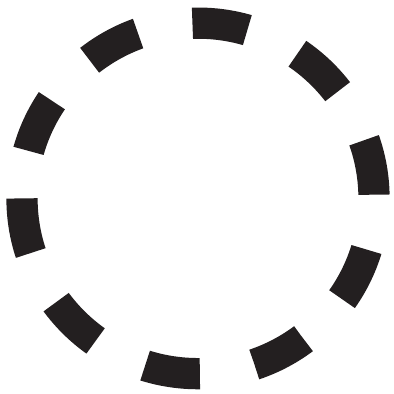}}
\newcommand{\ovdort}{\includegraphics [width=0.28cm]{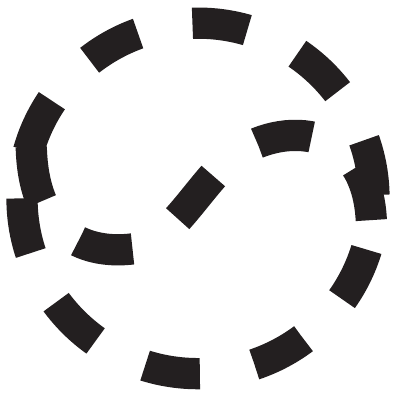}}

\newcommand{\ov}{\includegraphics [trim=0 0 0 0, clip, width=0.24cm]{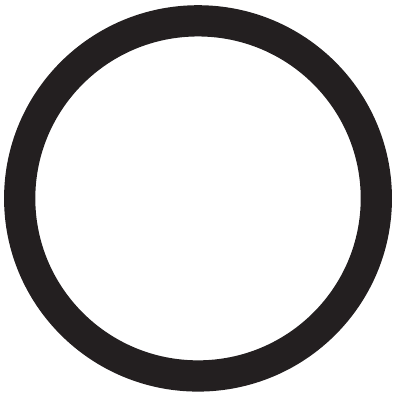}}
\newcommand{\ovkucuk}{\includegraphics [width=0.18cm]{o.pdf}}
\newcommand{\karkucuk}{\includegraphics [width=0.18cm]{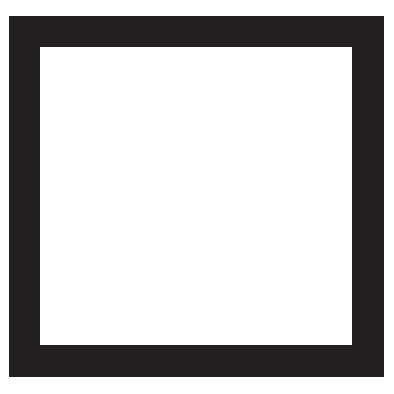}}
\newcommand{\kar}{\includegraphics [width=0.24cm]{k.pdf}}

\begin{document}

\title{Classification of totally real elliptic Lefschetz fibrations via necklace diagrams}

\author{Nerm\.{\i}n Salepc\.{\i}}
\address{Institut Camille Jordan,
Universit\'e Lyon I,
43, Boulevard du 11 Novembre 1918
69622 Villeurbanne Cedex, France}
\email{salepci@math.univ-lyon1.fr}


\begin{abstract}
We show that  totally real elliptic Lefschetz fibrations that admit a real section are classified by their ``real loci"  which is nothing but an $S^1$-valued Morse function on the real part of the total space. We assign to each such real locus a certain combinatorial object that we call a  \emph{necklace diagram}. On the one hand,  each necklace diagram  corresponds to an isomorphism class of a  totally real elliptic Lefschetz fibration that admits a real section, and on the other hand, it  refers to a decomposition of the identity into a product of certain matrices in $PSL(2,\Z)$.  Using an algorithm to find such decompositions, we obtain an explicit list of necklace diagrams associated with certain classes of totally real elliptic Lefschetz fibrations.  Moreover, we introduce refinements of necklace diagrams and show that refined necklace diagrams determine uniquely the isomorphism classes of the totally  real elliptic Lefschetz fibrations which may not have a real section.  By means of necklace diagrams we observe some interesting phenomena underlying special feature of real fibrations.
\end{abstract}

\maketitle


\section{Introduction}
As is well known Lefschetz fibrations are projections from an oriented connected smooth 4-manifold onto an oriented  connected smooth surface such that there exist finitely many critical points around which one can choose complex charts  so that  the projection on these charts  is given by $(z_{1}, z_{2}) \to z_{1}^2+z_{2}^2$. Regular fibers  of Lefschetz fibrations are oriented closed smooth surfaces of genus $g$, while singular fibers have only nodes. In the present work, we consider only those fibrations whose fiber genus is 1.  We call such fibrations  \emph{elliptic Lefschetz fibrations}.  Without loss of generality, we assume that  each singular fiber contains only one node and that no fiber contains a self intersection -1 sphere. 

The objects of our interest are  \emph{real elliptic Lefschetz fibrations} over $S^2$. They are defined as elliptic Lefschetz fibrations whose total and base spaces have \emph{real structures} which are compatible with the fiber structure.  A \emph{real structure} on an oriented smooth 4-manifold is defined as an orientation preserving involution whose fixed point set (which is called the \emph{real part}) has dimension 2, if it is not empty. It is worth mentioning here that not every 4-manifold admits such an involution.  Examples of  4-manifolds which do not admit  real structures can be found in \cite{KK}.
Likewise, a \emph{real structure} on a smooth oriented surface is defined as an orientation reversing involution. Obviously every surface admits a real structure.  Besides, the classification of real structures on surfaces up to conjugation by an orientation preserving diffeomorphism is known.  There are two invariants that determine the conjugacy class of a real structure:  its type (separating/ non-separating) and the number of the components of  its real part. A real structure is called \emph{separating} if the real part divide the surface into two disjoint halves; otherwise, it is called  \emph{non-separating}.  Throughout  the present work, the real structure considered on the base space $S^2$ will be the one induced from the complex conjugation on $\CP^1$. We denote it by $\cj$.  By definition of real Lefschetz fibrations, fibers over the real part $S^1$ of $\cj$ inherit real structures from the real structure of the total space.
We call such fibers \emph{real fibers}.  Real elliptic Lefschetz fibrations have 3 types of real regular fibers  that are classified by the number of  real components that can be 0,1, 2.  Only the structure with 2 real components is separating on $T^2$.

For the sake of simplicity, most of the time we assume that the real part $S^1$ is oriented. Fibrations with such a feature are called \emph{directed}.  Moreover, we consider mainly  fibrations which admit a real section (a section which commutes with the real structures of the total and the base space). But the cases of non-directed fibrations as well as of fibrations without a real section are also covered. The only essential condition imposed on fibrations is that all the critical values are real. Fibrations with only real critical values are called \emph{totally real}.

Our main interest is the topological classification of totally real elliptic Lefschetz fibrations.  Two real Lefschetz fibrations will be considered \emph{isomorphic} if  they can be carried one to other via orientation preserving equivariant diffeomorphisms. 
Recall that the classification of elliptic Lefschetz fibrations over $S^2$ has been known for over 30 years. It is due to Moishezon and Livn\'e \cite{Mo} that (non-real) elliptic Lefschetz fibrations over $S^2$ are classified by the number of critical values. The latter  is divisible by 12 and   
the class of elliptic Lefschetz fibrations with $12n$ critical values is denoted by $E(n)$, $n\in\N$.  
Furthermore,  $E(1)$ is isomorphic to the fibration $\CP^2 \#9 \overline{\CP^2} \to \CP^1$, obtained by blowing up a pencil of cubics in $\CP^2$, and 
 $E(n)= E(n-1)\sharp_{F} E(1)$ where $\sharp_{F}$ stands for the fiber sum of two fibrations.

In this note, we give the real version of this result for totally real elliptic Lefschetz fibrations.
The classification is obtained by means of certain combinatorial objects that we call \emph{(refined) necklace diagrams}. 
To each (refined) necklace diagram, we assign a \emph{monodromy}, a product of certain matrices in $PSL(2,\Z)$.  Indeed the product is always identity for fibrations over $S^2$, and what is crucial  is the decomposition of the identity. Necklace diagrams are  combinatorial counterparts of \emph{real Lefschetz chains} introduced in \cite{ner2}.   
Main results of this work, which are presented as Theorem~\ref{birebirkolye} and Theorem~\ref{birebirincikolye} covering the cases of directed totally real fibrations that admit a real section and  respectively fibrations possibly without  a real section, rely substantially on the material presented in \cite{ner2}.    As immediate corollaries of these theorems, we obtain that non-directed totally real elliptic Lefschetz fibrations admitting a section are classified by their necklace diagrams (defined up to symmetry and with the identity monodromy), while those fibrations which do not admit a real section are classified by the symmetry classes of refined necklace diagrams with the identity monodromy.  As a consequence of Theorem~\ref{birebirkolye},  we obtain an explicit list of totally real $E(1)$ and real $E(2)$ that admit a real section.  We investigate the algebraicity of  these fibrations and find the list of all real algebraic $E(1)$.  We also consider certain operations, \emph{mild/harsh sums}, \emph{flip-flops} and \emph{metamorphoses}, on the set of necklace diagrams. These operations allow us to construct new necklace diagrams from the given ones. By means of these operations, we construct examples of real Lefschetz fibrations which can not be written as the fiber sum of two real fibrations.


\textbf{Acknowledgements.} The material presented here is extracted from my thesis.  I am deeply indebted to my supervisors  Sergey Finashin and Viatcheslav Kharlamov  for  their guidance and limitless support.  I owe many thanks to Andy Wand who wrote the program to get the explicit list of necklace diagrams,  who also edited my present and past articles as a native english speaker.  
I  thank  Alex Degtyarev for his precious comments on the first manuscript and for productive discussions. 

The article has been finalized during my visit  to the Mathematisches Forschungsinstitut Oberwolfach as an ``Oberwolfach Leibniz fellow".
I would like to thank the institute for providing me exquisite working conditions.
 

\section{Real loci of real elliptic Lefschetz fibrations and necklace diagrams}

Let $\pi:X \to S^2$ be a directed real elliptic Lefschetz fibration. We look at the restriction, $\pi_{\R} : X_{\R}\to S^1$,  of $\pi$ to the real part $X_{\R}$ of $X$.  By definition,  fibers of  $\pi_{\R}$ are the real parts of the real fibers of $\pi$.  The base space $S^1$ is oriented (since we consider directed fibrations), whereas the total space $X_{\R}$  is either an empty set or a surface not necessarily oriented nor connected. 

By definition of real Lefschetz fibrations,  the map $\pi_{\R}$  is an  $S^1$-valued Morse function on $X_{\R}$ whose
regular fibers can be  $S^1$,  $S^1\amalg S^1$ or the empty set.  On the other hand,  singulars fibers  are either a wedge of two circles (this occurs in the case when the critical point is of index 1) or disjoint union of $S^1$ with an isolated point or  just an isolated point (these cases occur when the critical point is of index 0 or 2).  
As an immediate consequence, we note that the real part $X_{\R}$ is not empty if there is a real critical value.  
We consider only the fibrations with real critical values,  so $X_{\R}$ will never be empty throughout this article.

\begin{figure}[ht]
   \begin{center}
         \includegraphics[trim = 0 35 0 0, clip, width=3.6cm]{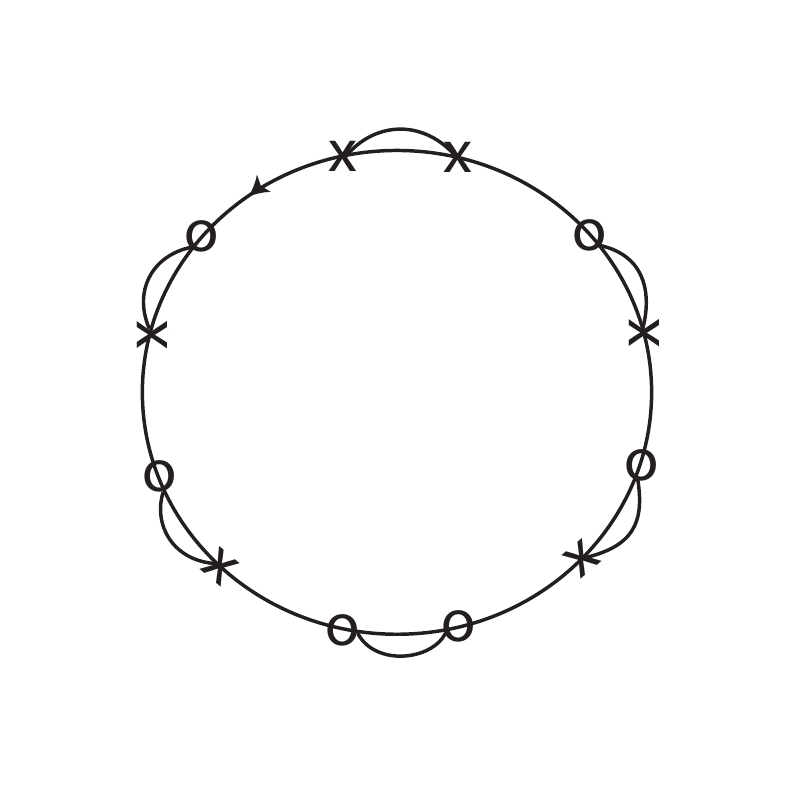}
      \caption{\small{Uncoated necklace diagram.}}
       \label{bebek}
       \end{center}
\end{figure}

For the sake of simplicity, we first focus on fibrations which admit a real section. By a real section, we understand a section $s:S^2\to X$ such that $c_{X}\circ s=s\circ \cj$.  
Now, let  $\pi_{\R} : X_{\R}\to S^1$ be the real locus of the directed real elliptic Lefschetz fibration $\pi:X\to S^2$ admitting a real section and having real critical values.  
We introduce a decoration on the base space $S^1$ as follows. First, we label the critical values of $\pi_{\R}$  by  ``$\times$"  or ``$\circ$"  according to the parity of indices of the corresponding critical points.  Namely, if the corresponding critical point is of index 1, we label the critical value by ``$\times$",  otherwise by  ``$\circ$".  (Note that $\pi_{\R}$ has critical values as long as $\pi$ has real critical values.)  
We now  consider a labeling on the set of  \emph{regular intervals},  $S^1\setminus \{\mbox{critical values}\}$, of $S^1$. Existence of a real section assures that fibers of   $\pi_{\R}$ are never empty, so there are only two possible topological types for regular fibers: $S^1$ or $S^1\amalg S^1$.  
Over each regular interval the topology of the fibers of $\pi_{\R}$ is fixed; moreover, it alternates as we pass through a critical value. 
We label regular intervals over which fibers have two components by doubling the interval, see Figure~\ref{bebek}. Regular intervals over which the fibers of $\pi_{\R}$ are a copy of $ S^1$ remain unlabeled.

\begin{figure}[ht]
   \begin{center}
         \includegraphics[trim = 5 5 0 0, clip, width=5cm]{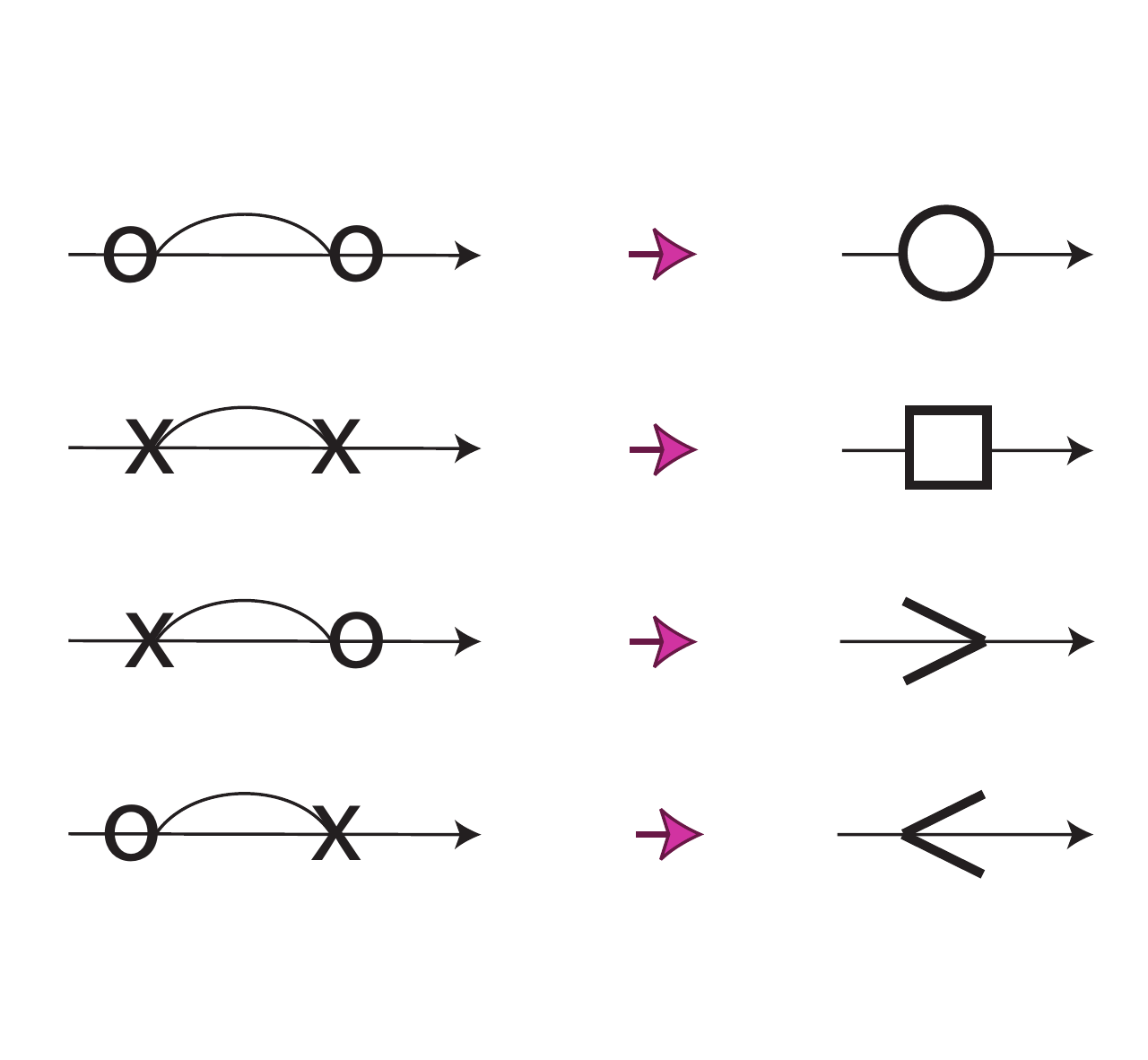}
      \caption{\small{Necklace stones}}
       \label{taslar}
       \end{center}
\end{figure}

Oriented $S^1$ together with such a decoration, is called an \emph{oriented uncoated necklace diagram}.
Let us now consider  ``standard" pieces of the uncoated necklace diagrams out of which we can built all possible uncoated necklace diagrams. To avoid the matching problem of  real structures, we deal with pieces of two critical values.  Let us choose a regular value on $S^1$.  (For some later use we choose the point on an unlabeled regular interval.)  With respect to this point and the orientation of $S^1$, we have 4 instances for a pair of two critical values. In order to simplify the decoration, for each instance we introduce new notations as shown in Figure~\ref{taslar}.

\begin{figure}[ht]
   \begin{center}
         \includegraphics[trim = 0 35 0 0, clip, width=10cm]{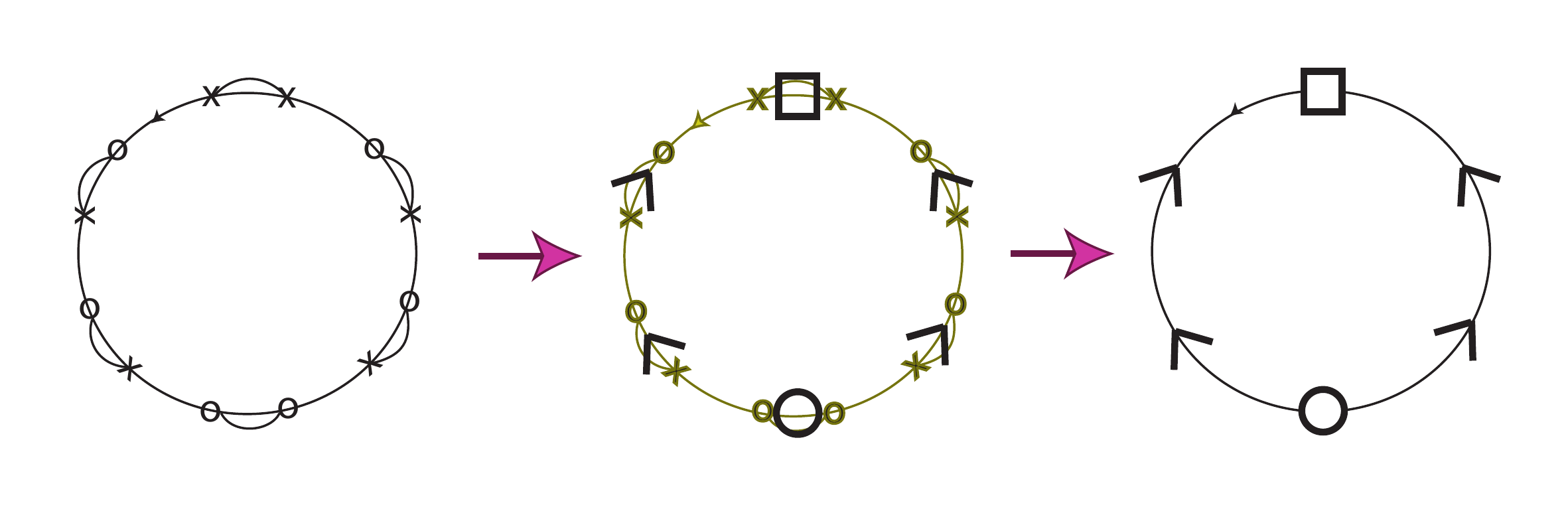}
      \caption{\small{A necklace diagram.}}
       \label{bebek2}
       \end{center}
\end{figure}

The oriented $S^1$ decorated using elements of the set  $\{\ov, \kar,>,<\}$ is called an \emph{oriented necklace diagram} (an example is shown in Figure~\ref{bebek2}). We call the elements of the set
$\{\ov, \kar,>,<\}$  \emph{(necklace) stones} and the pieces of the circle between the stones \emph{(necklace)
chain}.  Two oriented necklace diagrams are considered identical if they contain  the same types of stones going in the same cyclic order.

\begin{r.}
It is obvious from the construction that oriented necklace diagrams are invariants of directed real elliptic Lefschetz fibrations. 
When we consider non-directed real Lefschetz fibrations, we do not have a preferable orientation on the necklace diagram.  
Non-directed fibrations, hence, determine a pair of oriented necklace diagrams related by a mirror symmetry in which $\ov$-type and $\kar$-type stones remain unchanged, while $>$-type and  $<$-type stones interchanged.
\end{r.}


\section{Monodromy representations of stones}

As discussed in \cite{ner1},  monodromies of real Lefschetz fibrations around certain loops (namely, loops on which the real structure acts as a reflection)  break up into a composition of two real structures.  In particular,  the monodromy around a single real singular fiber can be written as $t_a=c'\circ c$  where $t_a$ denotes the positive Dehn twist along the vanishing cycle $a$ on a nearby marked (non-real) fiber $F$ identified with an abstract surface $T^2$,  and  $c', c:T^2\to T^2$ are the real structures  pulled from  the nearby right and respectively left real fibers (see \cite[Theorem~2] {ner1}).  Using this decomposition, to each decoration around a critical value, we assign a certain  $2\times 2$ matrix. These matrices are closely related to the monodromies of the fibrations.

 Let us first recall that each real structure,  $c:\Sa\to \Sa$, induces an isomorphism $c_*$ on $H_1(\Sa, \Z)=\Z \oplus \Z$ that  defines two rank 1 subgroups $H^c_{\pm}=\{\gamma \in H_{1}(T^2, \Z):c_*\gamma=\pm \gamma\}$ of $H_1(\Sa,\Z)$.  Moreover, if the real structure is separating, then $H_1(\Sa, \Z)=H^c_+ \oplus H^c_-$.  If $c$  is non-separating with one real component,  then $H_1(\Sa, \Z)/ H^c_+ \oplus H^c_-= \Z_{2}$.

\begin{r.}\label{vanissaykil}
By definition of real Lefschetz fibrations,  around each critical point and critical value we have equivariant local (closed) charts
$(U, \phi_U)$, $(V, \phi_V)$  such that $\pi|_U:(U,c_U)\to (V,\cj)$ is equivariantly isomorphic to either of $\xi_\pm:(E_\pm,\cj)\to (D_\epsilon,\cj)$, where 
$E_\pm=\{(z_1, z_2)\in \C^2\, :\, \left|z_1\right|\leq\sqrt{\epsilon},\,\left|z_1^2\pm z_2^2\right|\leq\epsilon^2\}$ and
$D_\epsilon=\{t\in \C\, :\, \left|t\right|\leq\epsilon^2 \},\, 0<\epsilon < 1$ with $\, \xi_\pm(z_1, z_2)=z_1^2\pm z_2^2$.
In the case of $\xi_{+}$ (this is the model for the critical point of index 0, 2) there are two types of real regular fibers distinguished by their real parts. 
In both cases,  one can choose invariant representatives for vanishing cycles and 
the action of the real structure  on the invariant representative can be either the antipodal map or the identity.  
On the other hand, in the case of $\xi_{-}$ (this is the case of critical points of index 1) real structure acts on the 
invariant representative of the vanishing cycle as a reflection.
Consequently, in the former situation ( which corresponds to the decoration ``$\circ$") the class of the vanishing cycle gives an element in $H^c_{+}$, while in the latter case (the case corresponding to the decoration ``$\times$") the class of the vanishing cycle gives an element of $H^c_{-}$. (A detailed discussion about the local models can be found in  of \cite[Section~3]{ner2} where the above claims are depicted in Figure~2.)  
\end{r.}

\begin{l.}\label{4P} For each decoration around a critical value,  we obtain the following matrices defined up to sign.
$$\begin{array}{llll}
P_{(-\times<)}&=\frac{1}{2}\small{\left(\begin{array}{cc}
           1 & 0\\
          -1 & 2
         \end{array} \right)},&P_{(>\times-)}&=\small{\left(\begin{array}{cc}
           2 & 0\\
          -1 & 1
         \end{array} \right)},\\

P_{(-\circ<)}&=\frac{1}{2}\small{\left(\begin{array}{cc}
           2 & 1\\
           0 & 1
         \end{array} \right)},& P_{(>\circ-)}&=\small{\left(\begin{array}{cc}
           1 & 1\\
           0 & 2
         \end{array}\right)}.
\end{array}$$
\end{l.}

\noindent {\it Proof: } The explicit calculations will be made for $P_{(-\times<)}$, the other cases are similar.

Let $q$ be the critical value decorated by ``$\times$". We consider a sufficiently small  $\epsilon$-neighborhood $D_{q}\subset S^2$ of  $q$ such that $\cj(D_{q})= D_{q}$ and $D_{q} \cap \{\mbox{critical values}\}=q$. We mark a non-real point $m$ on $\partial D^2$ and consider the shortest paths on $\partial D^2$ from $m$ to the left $q-\epsilon$ and right $q+\epsilon$ real points of $\partial D_{q}$.
By means of these paths, we pull the real structures on $F_{q\pm \epsilon}$ back to $F_{m}$.  Let us also fix an auxiliary identification $T^2 $ with the fiber $F_{m}$. 
Let $c:T^2\to T^2$ (respectively $c':T^2\to T^2$) be the real structures obtained by pulling back the real structure on  $F_{q- \epsilon}$ (respectively $F_{q+ \epsilon}$).

As discussed in Remark~\ref{vanissaykil},  the critical value of the type ``$\times$" provides a generator of $H_-^c\subset H_{1}(T^2,\Z)$. 
Let $b$ denote the corresponding vanishing cycle and $\beta$ the homology class of $b$.
We have $<\beta>=H_-^c$.  We choose a generator $\alpha$  for $H_+^c$ such that  $\alpha \circ \beta >0$.  As the decoration asserted, $c$ has one real component, so we want $\alpha \circ \beta=2$.

From the local monodromy decomposition, $t_b=c' \circ c$, we get  $c'_*={t_b}_*\circ c_*$; therefore,  $c'_*(\alpha)={t_b}_*( c_*(\alpha))=\alpha-2\beta$ and $c'_*(\beta)={t_b}_*( c_*(\beta))= -\beta$.   
Obviously,  the class $\alpha+ c'_*(\alpha)=2\alpha-2\beta  \in H_{+}^{c'} $, while $\beta-c'_*(\beta) \in H_{-}^{c'}$. 
We set  $\alpha'=\frac{\alpha-\beta}{2}$ and $ \beta'=\beta$ so that   ${\it{B'}}=\{\alpha', \beta'\}$ generates  $H_{+}^{c'}\oplus H_{+}^{c'}$.
The matrix $P_{-\times<}$ associated to the decoration ${-\times<}$ is, then,  $\frac{1}{2}\tiny{\left(\begin{array}{cc}
           1 & 0\\
          -1 & 2
         \end{array} \right)}$ which is the transition matrix from the base ${\it{B}}$ to the base ${\it{B'}}$.
         \hfill  $\Box$ \\

\begin{r.}\label{monmon}
Intuitively, the monodromy  assigned to a decoration around a critical value can be interpreted as the \emph{half} of the monodromy around the real critical value, see  Figure~\ref{yarim}. 
We observe that  $P_{>\times-}P_{-\times<}= [t_{b}]_{\it{B}}$ where $[t_{b}]_{\it{B}}$ denotes the matrix of $t_{b*}$ with respect to the base ${\it{B}}$. Similarly,  we have $P_{>\circ-}P_{-\circ<}= [t_{a}]_{\it{B}}$ as well as  $P_{-\times<}P_{>\times-}= [t_{b}]_{\it{B'}}$ and $P_{-\circ<}P_{>\circ-}= [t_{a}]_{\it{B'}}$.  Moreover, we have  $P_{-\times<}= M P^{-1}_{>\times-} M$ and  $P_{-\circ<}=M P^{-1}_{>\circ-} M$ where $M=M^{-1}={\tiny{\left(\begin{array}{cc}
           1 & 0\\
           0 & -1
         \end{array} \right)}}$.  
         
         \begin{figure}[ht]
\begin{center}
      \includegraphics[trim = 10 30 10 10, clip, width=2cm]{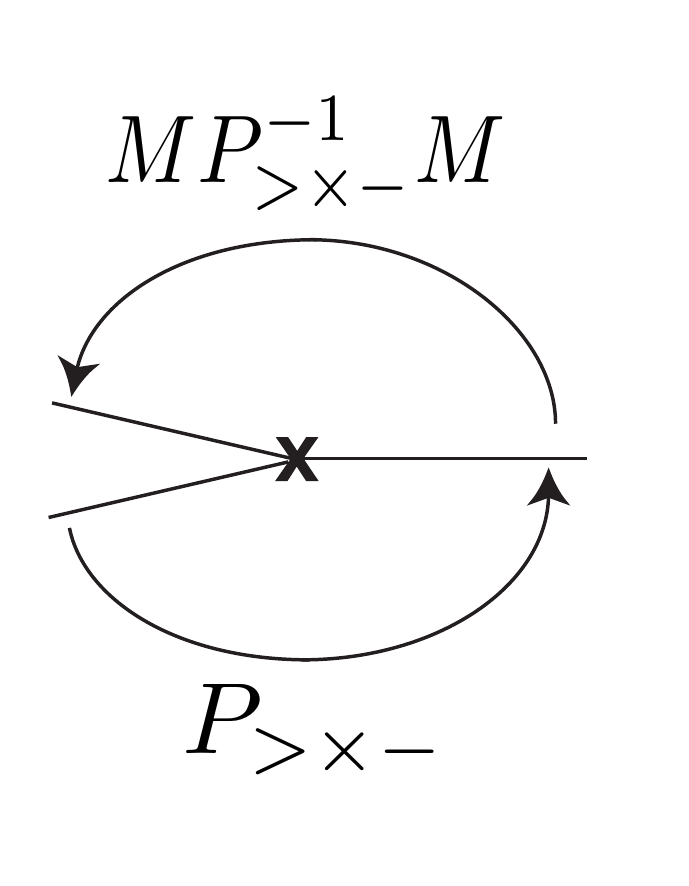}
       \caption{\small{Decomposition of the monodromy associated to the decoration of a real critical value.}}
       \label{yarim}
       \end{center}
 \end{figure}
         
\end{r.}

Therefore, to each necklace stone, we assign the following products (defined up to sign). 

$$\begin{array}{lll}
P_{\karkucuk}&= P_{(-\times<)}P_{(>\times -)}&=\small{\left(\begin{array}{cc}
           1 & 0\\
           -2 & 1
         \end{array} \right)},\\
         &&\\
P_{\ovkucuk}&=P_{(-\circ<)} P_{(>\circ-)}&=\small{\left(\begin{array}{cc}
           1 & 2\\
           0 & 1
         \end{array} \right)},\\
         &&\\
P_{>}&=P_{(-\times<)} P_{(>\circ-)}&=\small{\frac{1}{2}\left(\begin{array}{cc}
           1 & 1\\
           -1 & 3
         \end{array} \right)},\\
         && \\
P_{<}&=P_{(-\circ<)} P_{(>\times-)}&=\small{\frac{1}{2}\left(\begin{array}{cc}
           3 & 1\\
           -1 & 1
         \end{array} \right)}.
\end{array}$$

It is known that mapping class group, $\M(T^2)$ is isomorphic to the group $SL(2,\Z)$. 
Let $T^2$ be identified with $S^1\times S^1$.  We set  $a=S^1\times \{0\}$ and $b=\{0\} \times S^1$. 
Then, we consider two presentations of $\M(T^2)\simeq SL(2,\Z)$ as follows. 

$$\begin{array}{lcl}
SL(2,\Z)&=&\{[t_{a}]={\tiny{\left(\begin{array}{cc}
           1 & 1\\
           0 & 1
         \end{array} \right)}} \,\,\textrm{and}\,\,[t_{b}]=\tiny{\left(\begin{array}{cc}
           1 & 0\\
           -1 & 1
         \end{array} \right)} : ([t_{a}][t_{b}])^6=id\}\\
         
        &=& \{x={\tiny{\left(\begin{array}{cc}
           0 & 1\\
          -1 & 0
         \end{array} \right)}} \,\,\textrm{and}\,\, y=\tiny{\left(\begin{array}{cc}
           0 & 1\\
          -1 & 1
         \end{array} \right)} : x^2=y^3, x^4=id\}. \end{array}$$ 
         
Let us note one can switch from the first presentation to the second by setting $x= [t_{a}] [t_{b}] [t_{a}]=[t_{b}] [t_{a}] [t_{b}]$ and $y=[t_{a}] [t_{b}]$.
Since $x^2= -id$, we have $PSL(2,\Z)=\{x, y: x^2=y^3=id\}$.

\begin{r.}
Because we choose  a point on an unlabeled interval, the matrices have coefficients in $\frac12 \Z$.  If we marked a regular value on an labeled interval, the matrices we get would be elements of $PSL(2,\Z)$.  
(The reason why we prefer a point on an unlabeled interval is to get a nice relation between necklace stones and the real part (see  Remark~\ref{tophop}). )
The subgroup generated by $\{P_{\ovkucuk}, P_{\karkucuk}, P_{>}, P_{<}\}$  is conjugate to $PSL(2, \Z)$. 
To be able work with $PSL(2,\Z)$, we consider the following lemma.
\end{r.}

\begin{l.} Let $R=\frac{1}{2}\tiny{\left(\begin{array}{cc}
           1 & -1\\
           1 & 1
         \end{array} \right)}$ and $\mathbb{P}=R^{-1}PR$.
         
         Then, for each necklace stone we obtain the following factorization.

$\begin{array}{ll} \label{4Ptilde}
{\mathbb{P}}_{\karkucuk}&= yxy\\
{\mathbb{P}}_{\ovkucuk}&= xyxyx\\
{\mathbb{P}}_{>}&=y^2x\\
{\mathbb{P}}_{<}&=xy^2
\end{array}$
\end{l.}

\noindent {\it Proof: }
The proof follows  from the observation  that  
${\mathbb{P}}_>={[t_{a}]},\,\,{\mathbb{P}}_<={[t_{b}]}$, while ${\mathbb{P}}_{\karkucuk}={[t_{a}] [t_{b}] [t_{a}]^{-1}},\,\,{\mathbb{P}}_{\ovkucuk}={[t_{a}]^{-1}[t_{b}] [t_{a}]}.$ 
 \hfill  $\Box$ \\

The matrices  $\mathbb{P}_{\ovkucuk}, \mathbb{P}_{\karkucuk}, \mathbb{P}_{>}, \mathbb{P}_{<}$ are called  \emph{monodromies of stones}. The product of monodromies of necklace stones is called the \emph{monodromy of a necklace diagram}.

\begin{l.}\label{monod}
Let $\pi: X\to S^2$ be a directed totally real elliptic Lefschetz fibration admitting  a real section. Then, the monodromy of the oriented necklace diagram associated with
$\pi$ is the identity in $PSL(2,\Z)$.
\end{l.}

\begin{figure}[ht]
\begin{center}
      \includegraphics[scale=0.3]{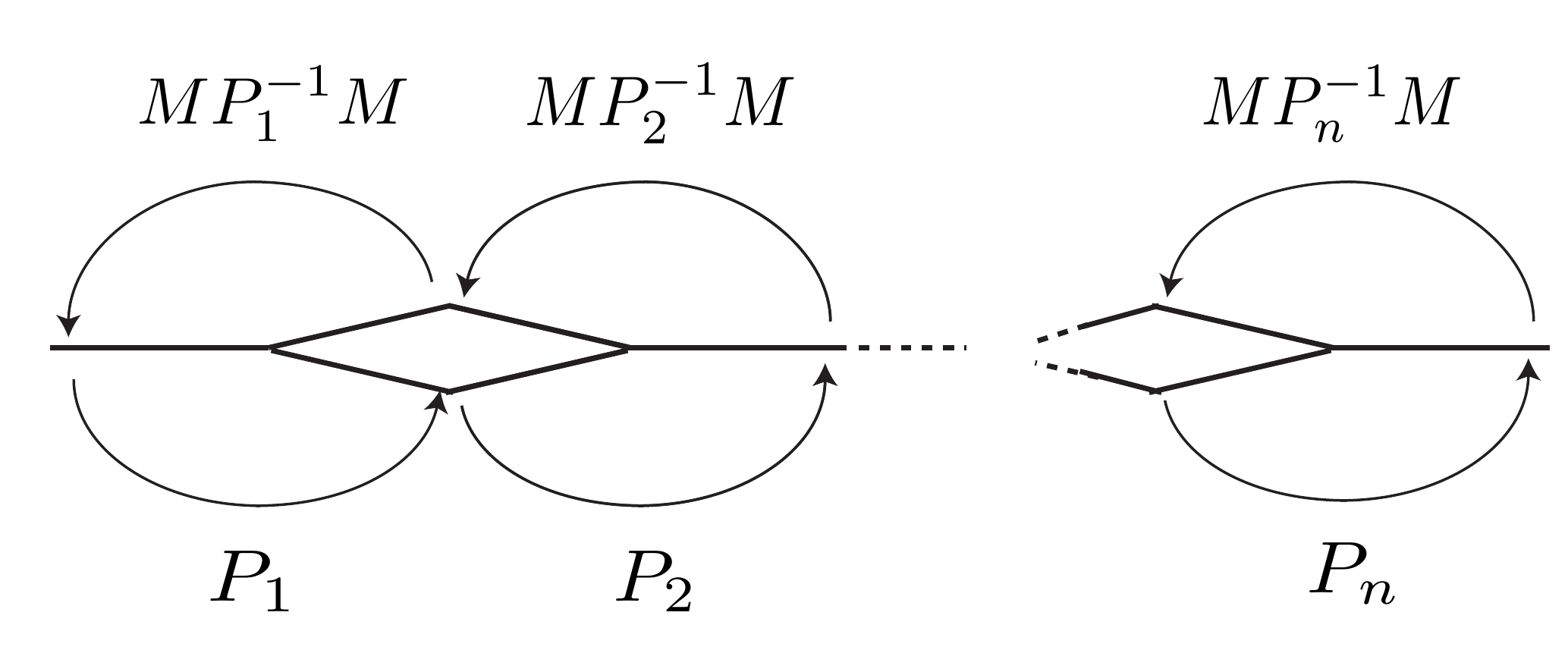}
       \caption{\small{Total monodromy as the product  $P_{1}.P_{2},\ldots P_{n}$.}}
       \label{tamtam}
       \end{center}
 \end{figure}

\noindent {\it Proof: } 
Let $\{q_{1},\ldots, q_{n}\}$ be the ordered
set of real critical values of $\pi$. By means of underlying uncoated necklace diagram, we can write the monodromy 
of the necklace diagram as the product  $P=P_{1}P_{2}\ldots P_{n}$ where $P_{i}$
is  the matrix associated to the decoration around  $q_{i}$. Following Remark~\ref{monmon},  we can write  the monodromy  along a curve surrounding all real critical values can be written as the product (see Figure~\ref{tamtam}) $$P_{1}P_{2}\ldots P_{n}{\tiny{\left(\begin{array}{cc}
           1 & 0\\
           0 & -1
         \end{array} \right)}}P_{n}^{-1}\ldots P_{2}^{-1} P_{1}^{-1} {\tiny{\left(\begin{array}{cc}
           1 & 0\\
           0 & -1
         \end{array} \right)}}.$$ 
         
          If there is no non-real critical value, the monodromy along the curve we consider is identical to the total monodromy  of the fibration which is the identity.
          
         Thus, we get  $P{\tiny{\left(\begin{array}{cc}
           1 & 0\\
           0 & -1
         \end{array} \right)}}P^{-1} {\tiny{\left(\begin{array}{cc}
           1 & 0\\
           0 & -1
         \end{array} \right)}}=id$.  The equality assures that $P$ is the identity in $PSL(2,\Z)$.
 \hfill  $\Box$ \\

\begin{r.} \label{dual}
An important observation is that ${\mathbb{P}}_{\ovkucuk}=x{\mathbb{P}}_{\karkucuk} x$ and
${\mathbb{P}}_<=x {\mathbb{P}}_{>} x$. Hence, if a necklace diagram has
the identity monodromy, then the necklace diagram obtained from the original by replacing
each $\kar$-type stone with $\ov$-type stone, and each
$>$-type stone with $<$-type stone and vice versa, has also the identity monodromy. Necklace diagrams obtained in this manner are called \emph{dual necklace diagrams}.
 \end{r.}

\begin{r.}  \label{mono} Although for the moment,  we  focus on fibrations admitting a real section, it is worth mentioning the case of real structures with no real component.  Indeed, real structures with no real component and those structures with two real components are are isotopic to each other as orientation reversing diffeomorphisms, although they are two non-isotopic real structures. As a result, they induce the same isomorphism on the homology group; hence, monodromy calculations remain the same if  the real structure with two real components is replaced by a real structure with no real component.
\end{r.}


\section{The correspondence theorems and consequences}

\begin{t.} \label{birebirkolye}
There exists a one-to-one correspondence between the set of oriented necklace diagrams with $6n$ stones whose monodromy is the
identity and the set of isomorphism classes of directed totally real elliptic Lefschetz fibrations $E(n)$,  $n\in\N$, that admit a real section.
\end{t.}

\noindent {\it Proof: }
From the discussions  of  the previous sections and by  Lemma~\ref{monod}, we have an injection the set of fibrations to the set of diagrams.
We want to show that this map is well-defined and surjective.  

The crucial observation is that each necklace diagram with identity monodromy defines uniquely, up to cyclic order, a  \emph{real Lefschetz chain of conjugacy classes of real codes}. Let us recall that a \emph{real code} is  a pair $(c, a)$ consisting of a simple closed curve $a$  and a real structure $c$ on the fiber such that $c(a)=a$.  As shown in \cite{ner2}, conjugacy classes of real codes are complete invariants of  equivariant neighborhoods of real singular fibers of a real Lefschetz fibration.   A \emph{real Lefschetz chain} of conjugacy classes of real codes is a chain of codes  $\{c_{1}, a_{1}\}, \{c_{2}, a_{2}\}, \ldots\{c_{n}, a_{n}\}$ such that $c_{i+1}$ is conjugate to $t_{a_{i}}\circ c_{i}$.   

To understand the relation between a necklace diagram and a real Lefschetz chains, it is enough to investigate the decorations of  the underlying uncoated necklace diagram. 
By definition, the decoration on regular intervals determine the conjugacy classes of real structures on the fibers over this interval.
By Remark~\ref{vanissaykil},  the decoration on the critical value determines the isotopy class of the vanishing cycle (invariant under the action of the real structure). The decoration around a critical value, hence, dictates the conjugacy class of  a real code.  An oriented necklace diagrams, hence, defines an ordered sequence of the conjugacy classes of real codes, so defines a real Lefschetz chain up to cyclic order. The result,  therefore, follows from \cite[Proposition~16]{ner2}.
 \hfill  $\Box$ \\

\begin{r.} As mentioned in  Remark~3.6 of \cite{ner2}, in the case when the regular fibers are tori, there are 6 conjugacy classes of real codes.
 Whereas, around each critical value we have 4 different decorations $\{-\times<, >\times-, -\circ<, >\circ -\}$.  These decorations are exactly 4 of the 6 possible real codes on $T^2$. 
The other two cases appear in the case of non-existence of a real section that we discuss in the last section. 
\end{r.}

\begin{c.} \label{birebirsimlikolye}
There exists a bijection between the set of symmetry classes of non-oriented necklace diagrams with $6n$ stones whose
monodromy is the identity, and the set of isomorphism classes of
non-directed totally real $E(n)$, $n\in\N$ which admit a real section.
$\Box$
\end{c.}

Each necklace diagram defines a decomposition of the identity in $PSL(2,\Z)$ into a product of $6n$ elements that  are chosen from the set of monodromies of necklace stones.  
There is a simple algorithm to find all necklace diagrams associated with $E(n)$.  Applying the algorithm, we obtain the complete list of necklace diagrams of $E(1)$. 
Later Andy Wand wrote a computer program for  $n=1, 2$.

The following theorem concerns $n=1$. 

\begin{t.}
There exist precisely 25 isomorphism classes of 
non-directed  totally real $E(1)$ admitting a real section. These classes are characterized
by the non-oriented necklace diagrams presented in Figure~\ref{e1kolyeler}.
\begin{figure}[ht]
\begin{center}
      \includegraphics[trim = 0 0 0 0, clip, width=12cm]{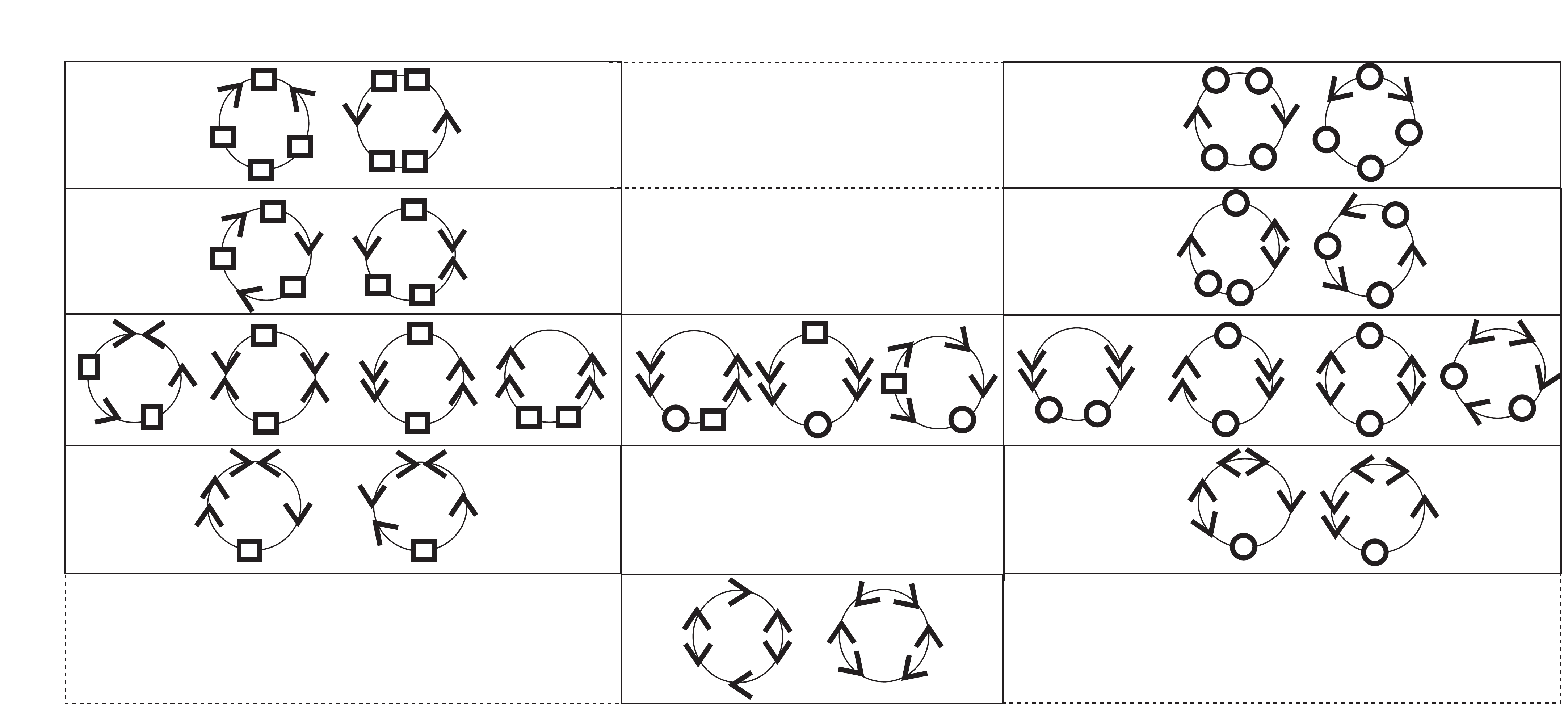}
       \caption{\small{List of necklace diagrams of totally real $E(1)$ admitting a real section.}}
       \label{e1kolyeler}
       \end{center}
 \end{figure}
\end{t.}

\noindent {\it Proof: }
By Theorem~\ref{birebirsimlikolye}, it is enough to find the list of
symmetry classes of necklace diagrams of $6$ stones whose
monodromy is the identity.  Each necklace diagram defines a decomposition of the identity in $PSL(2,\Z)=\{x, y: x^2=y^3=id \}$ into a product of elements $yxy,  xyxyx,  xy^2, y^2x$.
Let   $S=yxy, C=xyxyx, L=xy^2, R=y^2x$. To find all decompositions, we first consider all the words of length $6$ of letters $S,C,L,R$ such that the product is the identity.
We, then, quotient out the words which are equivalent to each other up to cyclic ordering. This way we obtain all oriented necklace diagrams which have the identity monodromy. 
As a final step, we quotient out the symmetry classes.  In terms of words composed of the letters  $S, C, L, R$,  symmetry classes  can be interpreted as follows:  two words are considered to be equivalent if one is the reversed of the other with each $L$  is replaced by $R$, and vice versa.  For example, $CLLSRR \sim LLSRRC$.  
 \hfill  $\Box$ \\

It is worth mentioning here that there are 8421 many necklace diagrams of 12 stones with the identity monodromy.
Below, in Proposition~\ref{propmaxE2}, we give  explicit list of necklace diagrams corresponding to certain classes.  Later, we also explore some interesting examples.

\begin{r.} \label{tophop} The topological invariants of $X_{\R}$ can be read from the necklace diagram of  $\pi:X \to S^2$. 
Namely, we have $\beta_{0}(X_{\R})=\beta_{2}(X_{\R})= |\ov|+1$ and  $\beta_{1}(X_{\R})= 2(|\kar|+1)$ where $\beta_{i}$ denotes the $i^{th}$ Betti number of  $X_{\R}$ and 
$|\ov|$, $|\kar|$  denote the number of $\ov$-type and  respectively $\kar$-type stones of the necklace diagram associated with $\pi_{\R}$. 
Consequently,  we have the Euler characteristic  $\chi(X_\R)=2(|\ov|-|\kar|)$,  and the total Betti number $\beta_*(X_\R)=2(|\ov|+|\kar|)+4$. 
 \end{r.}

Recall that in general we have   $\beta_*(X_{\R})\leq \beta_*(X)$ (known as \emph{Smith inequality}).
It is known that  $\beta_*(E(n))=12n$ (\cite{GS}), so we have   $\beta_*(E(n)_{\R})\leq 12n$. 

\begin{d.} A real structure $c_{X}$ on $X$ is called \emph{maximal}  if $\beta_*(X_{\R})=\beta_*(X)$.
A necklace diagram \emph{maximal} if  $|\ov|+|\kar|=\frac{12n-4}{2}$. 
\end{d.}

We have the following immediate consequences.

\begin{p.}\label{4arrow} Each  necklace diagram whose monodromy is the identity contains at least two arrow type stones. \hfill  $\Box$ \\
\end{p.}

\begin{c.}\label{enazbirx}
A totally real elliptic Lefschetz fibration admitting a real section contains at least two critical values of type $``\times".$
$\Box$
\end{c.}

There are 4 maximal $E(1)$ whose necklace diagrams that  are depicted on the top line of Figure~\ref{e1kolyeler}.  For $n=2$ we have:

\begin{p.}\label{propmaxE2}
There  are 10 isomorphism classes of  maximal non-directed totally real $E(2)$ admitting a real section. 
Corresponding necklace diagrams are given Figure~\ref{maxE2}. \hfill  $\Box$ \\

\begin{figure}[ht]
\begin{center}
 \includegraphics[trim = 0 0 0 0, clip, width=9cm]{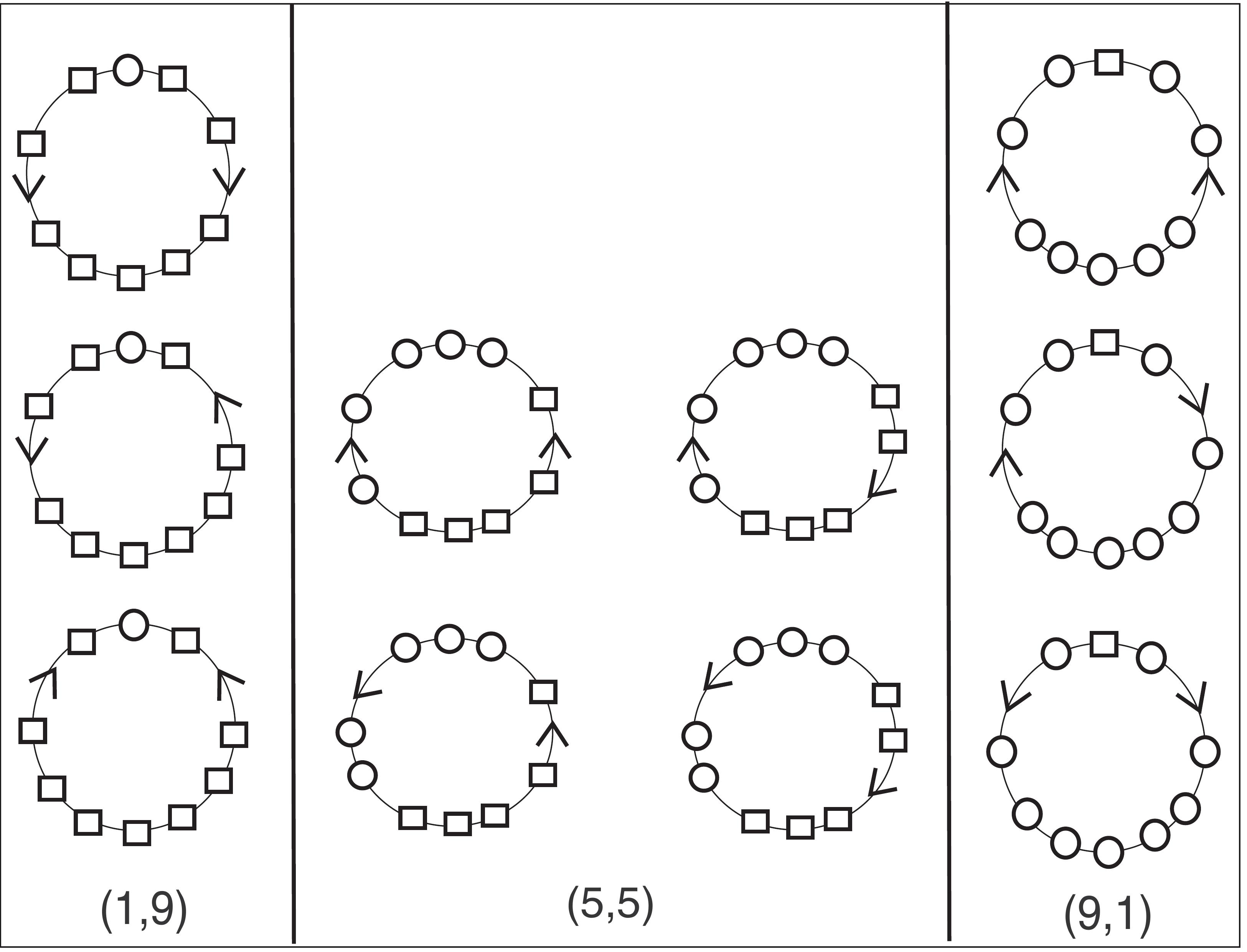}
       \caption{\small{List of necklace diagrams of maximal totally real $E(2)$ admitting a real section.}}
       \label{maxE2}
       \end{center}
 \end{figure}
\end{p.}


\section{Applications of necklace diagrams}

In this section, we study the algebraic realization of the totally real elliptic
Lefschetz fibrations admitting a real section. 
The crucial observation is that any algebraic elliptic Lefschetz fibration $E(n)$ admitting a real section 
can be seen as the double branched covering of a Hirzebruch surface of degree $2n$, branched
at the exceptional section and a trigonal curve disjoint from the
section.  Orevkov~\cite{orev} introduced a 
 real version of Grothendieck's dessins d`enfants for  the trigonal curves, which are
disjoint from the exceptional section, on Hirzebruch surfaces.  
We apply his result by  converting language of real dessin d'enfants to the language of necklace diagrams.

\subsection{Trigonal curves on Hirzebruch surfaces}

The Hirzebruch surface, $H(k)$, of
degree $k$ is a complex surface equipped with a projection,
$\pi_k:H(k)\to \CP^1$, which defines a $\CP^1$-bundle over $\CP^1$
with a unique exceptional section $s$ such that
$s\circ s= -k$. In particular, $H(0)=\CP^1 \times \CP^1$ and
$H(1)$ is $\CP^2$ blown up at one point.

Each Hirzebruch surface $H(k)$ can be obtained from $H(0)$ 
by a sequence of blow-ups followed by blow-downs at a certain set of points. If
these points are chosen to be real, then the resulting Hirzebruch
surface has a real structure inherited from the real structure
$\cj\times\cj$ on $H(0)$. This will be  the real structure of our consideration.
With respect to this real structure, the real part of $H(k)$ is a torus if $k$ is even; otherwise
it is a Klein bottle.

In this note, we only consider nonsingular curves, so by a
\emph{trigonal curve} on a Hirzebruch surface $H(k)$ we understand
a smooth algebraic curve $C\subset H(k)$ such that the restriction
of the bundle projection, $\pi_k:H(k)\to \CP^1$, to $C$ is of
degree 3. A trigonal curve on $H(k)$ is called \emph{real} if it
is invariant under the real structure of $H(k)$.

\subsection{Real dessins d'enfants of  trigonal curves}

We choose affine (complex) coordinates $(x, y)$ for $H(k)$ such that the
equation $x=\mbox{\emph{const}}$ corresponds to fibers of $\pi_k$ and
$y=\infty$ is the exceptional section $s$. Then, with respect to
such affine coordinates any (algebraic) trigonal curve
can be given by a polynomial of the form $y^3+u(x)y+v(x)$ where
$u$ and $v$ are real one variable polynomials such that $deg\,u=2k$
and $deg\,v=3k$.

The discriminant of $y^3+ u(x)y+v(x)=0$ with
respect to $y$ is $-4u^3-27v^2$.
Let  $D=4u^3+27v^2$. The fraction
$j=\frac{4u^3}{D}$ is the \emph{$j$-invariant} of a trigonal curve $C \subset H(k)$.
The $j$-invariant defines a real rational function $j:\CP^1\to \CP^1$ whose poles
are the roots of $D$,  zeros are the
roots of $u$ (taken with multiplicity 3),  and the solutions of $j=1$ are the roots of $v$
(taken with the multiplicity 2).

\begin{figure}[ht]
\begin{center}
\includegraphics[scale=0.3,trim=0 0 70 0]{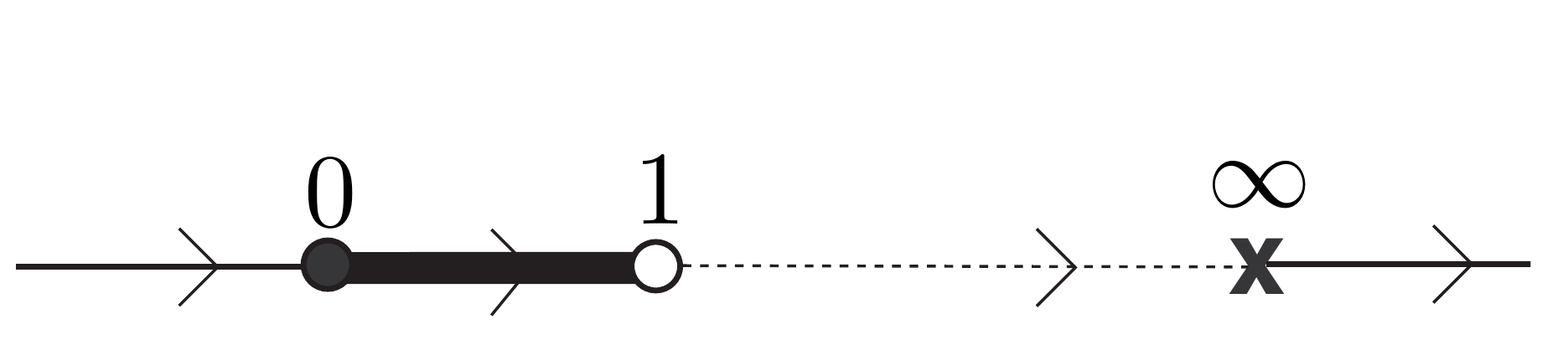}
\caption{\small{Coloring of $\RP^1$.}}
\label{ApRP1}
\end{center}
\end{figure}

\begin{figure}[ht]
\begin{center}
\includegraphics[scale=0.3,trim=0 0 70 0]{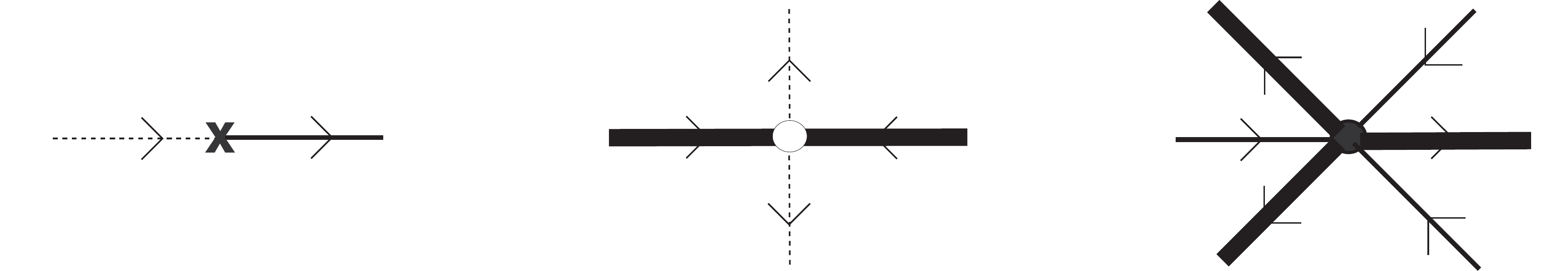}
\caption{\small{The graph around the inverse images of zeros of $ D, v, u$.}}
\label{Apkoseler}
\end{center}
\end{figure}

Let us color $\RP^1$ as in  Figure~\ref{ApRP1}.
Then, the inverse image $j^{-1}(\RP^1)$ turns naturally into
an oriented colored graph on $\CP^1$.  Since $j$ is real, the graph is symmetric with respect to the complex conjugation on $\CP^1$.
Around the vertices  the graph looks as shown Figure~\ref{Apkoseler}. (Detailed discussion on $j$-invariant of trigonal curves can be found in  cf.  \cite{DIK}, \cite{orev}.)

The following theorem gives the conditions which are sufficient for the (real) algebraic realizability of a graph and the existence of respective polynomials $u, v, D$.

\begin{t.}\cite{orev}\label{torev}
Let $\Gamma\subset S^2$ be an embedded oriented graph where  each of its edges is one of the three kinds: $\includegraphics [width=0.6cm]{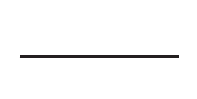}, \includegraphics [width=0.6cm]{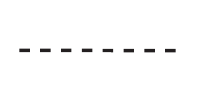}, \includegraphics[width=0.6cm]{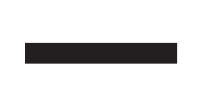}$ and some of its vertices are colored by the elements of the set $\{\circ, \bullet, \times\}$, while others remains uncolored. Let $\Gamma$ satisfy the following conditions:\\
(1) The graph $\Gamma$ is symmetric with respect to an equator of $S^2$, which is included
into $\Gamma$;\\
(2) The valency of each vertex $``\bullet"$ is divisible by $6$, and the incident edges are colored alternatively by incoming $\includegraphics [width=0.6cm]{Apkenar1.pdf}$, and outgoing $\includegraphics [width=0.6cm]{Apkenar3.pdf}$;\\
(3) The valency of each vertex $``\circ"$ is divisible by 4, and the incident edges are colored alternatively by incoming $\includegraphics [width=0.6cm]{Apkenar3.pdf}$, and outgoing $\includegraphics [width=0.6cm]{Apkenar2.pdf}$;\\
(4) The valency of each vertex $``\times"$ is 2, and the incident edges are colored alternatively by incoming $ \includegraphics [width=0.6cm]{Apkenar2.pdf}$, and outgoing $\includegraphics [width=0.6cm]{Apkenar1.pdf}$;\\
(5) The valency of each non-colored vertex is even, and the incident edges are of the same color;\\
(6) Each connected component of $S^2\setminus \Gamma$ is homeomorphic to an open disc whose boundary is colored as a
covering of $\RP^1$ (colored and oriented as in Figure~\ref{ApRP1}) and the orientations of the boundaries of neighboring discs are opposite.\\
Then, there exists a real rational function $j=\frac{4u^3}{D}$ whose graph is $\Gamma$. 
(And thus, there exist a non-singular real algebraic trigonal curve associated to the  $j$ invariant.)
\end{t.}

\begin{d.}
A graph on $S^2$ satisfying the conditions (1)-(6) of the above
theorem is called a \emph{real dessin d'enfant}.
\end{d.}

\begin{r.}
Let us accentuate the fact that there is no relation between the decorations ``$\times", ``\circ$" of the critical values that we use to introduce the necklace diagrams  and the coloring of the vertices of the real dessin d'enfants considered in this section.     
\end{r.}

\subsection{Correspondence between real schemes and real dessins d'enfants}

\begin{figure}[ht]
\begin{center}
        \includegraphics[scale=0.16,trim=0 0 00 0]{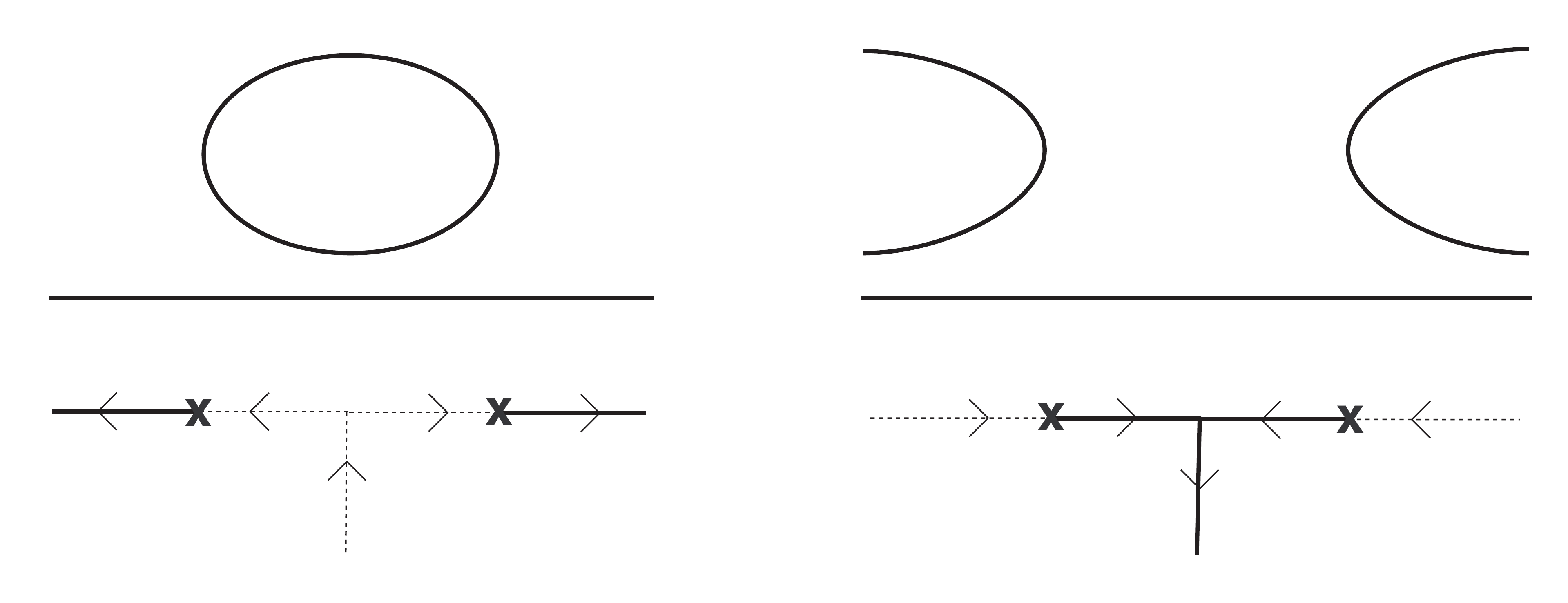}
        \includegraphics[scale=0.16,trim=0 0 00 0]{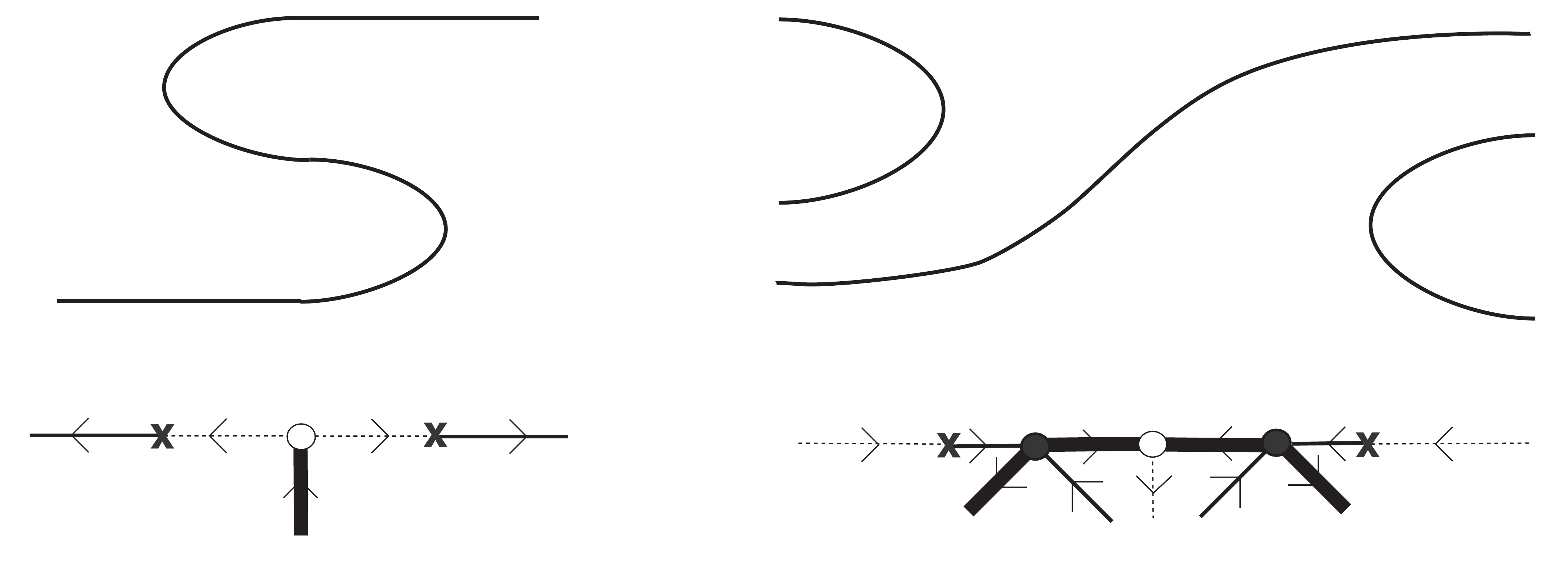}
\caption{\small{Segments of the curve corresponding to fragments of the minimal graphs.}}
\label{Aparcalar}
\end{center}
\end{figure}

The real scheme of a trigonal curve imposes strong restrictions on the arrangement of the real roots of $u$, $v$ and $D$.
For example, the zeros of $D$ correspond to the points where the trigonal curve is tangent to the fibers of $\pi_k :H(k)\to \CP^1$.
A typical correspondence for certain model pieces of the curve is shown in Figure~\ref{Aparcalar}.
Because the graph is symmetric with respect to the equator, we consider the part of the graph lying on one of the half discs. 

As we mention in the previous section, necklace diagrams encode the topology of the real part (except orientability) of  $E(n)$ which admit a real section.
Indeed, the real part of totally real $E(n)$,  admitting a real section,  consists of  spherical components (the number of which is $|\ov|$) and a higher genus component which is an orientable surface of genus  $|\kar|+1$  if $n$ is even; a non-orientable surface with $2|\kar|+1$ cross-caps, otherwise.

\begin{figure}[ht]
\begin{center}
\includegraphics[scale=0.32,trim=0 0 0 0, clip]{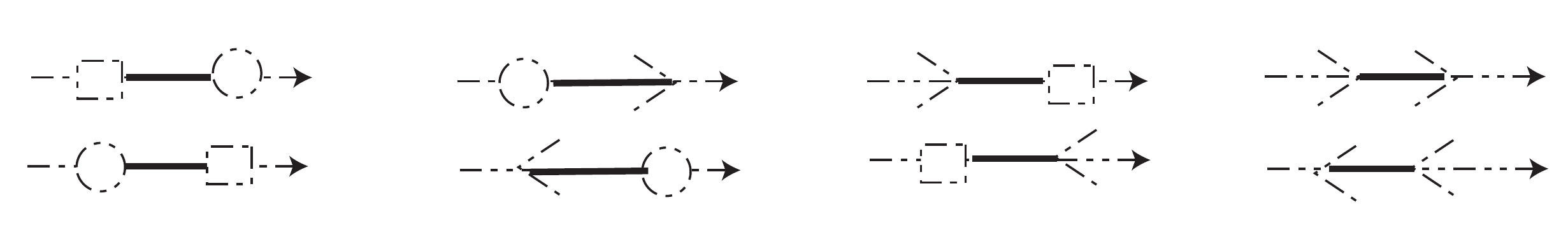}
\caption{\small{Essential intervals.}}
\label{Apkolyeilis}
\end{center}
\end{figure}

\begin{d.}
A segment of a necklace diagram is called \emph{essential} if the corresponding graph fragment  contains at least one ``$\circ$" type vertex and at least two ``$\bullet$" type vertices. (Essential segments are listed in Figure~\ref{Apkolyeilis}.) 
\end{d.}

\subsection{Applications}

\begin{p.}\label{yasak}
If a real elliptic Lefschetz fibration, $E(n)$, admitting a real section is algebraic then the corresponding necklace diagram has the following properties:

\begin{itemize}
\item there are not more than $2n$ essential segments,

 \item the sum of the number of essential segments and the number of arrow type stones cannot be greater then $6n$.
\end{itemize}
\end{p.}

\noindent {\it Proof: } 
For a trigonal curve on $H(2n)$ defined by $y^3+u(x)y+v(x)$,  $deg\, u=2\cdot 2n$ and $deg\, v= 3\cdot 2n$.
Thus, the real dessin d'enfant can have at most $4n$ vertices colored by ``$\bullet$'' and at most $6n$ vertices colored by ``$\circ$''.
The result follows from the observation that each essential interval corresponds to a graph fragment which contains
at least two ``$\bullet$'' type vertices and at least one ``$\circ$'' type vertex, while
each arrow type stone corresponds  to a fragment having  at least one ``$\circ$'' type vertex.
 \hfill  $\Box$ \\

\begin{c.}
The totally real elliptic Lefschetz fibrations corresponding to the necklace diagrams depicted in Figure~\ref{Apcebdegil} are not realized algebraically.
 \begin{figure}[ht]
\begin{center}
       \includegraphics[scale=0.33,trim=0 0 -40 0]{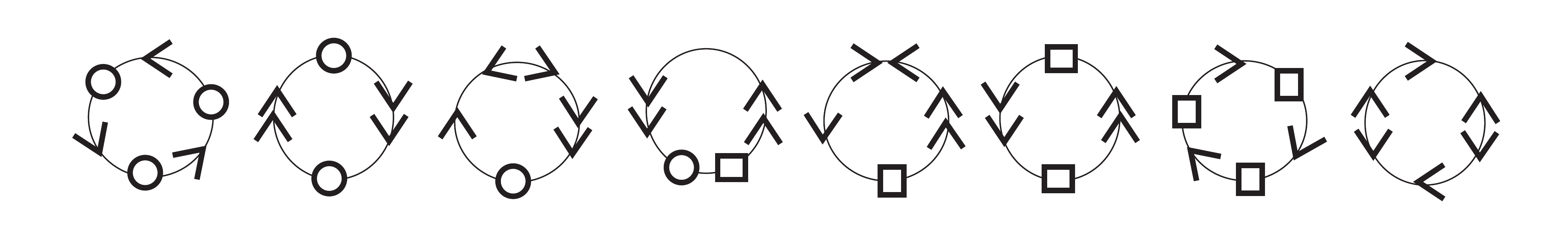}
\caption{\small{Necklace diagrams violating algebraicity.}}
\label{Apcebdegil}
\end{center}
\end{figure}
 \end{c.}

\noindent {\it Proof: }
It is easy to see that the diagram having only arrow type stones violates the second condition, while all the others violate the first condition stated in Proposition~\ref{yasak}.
\hfill  $\Box$ \\ 

\begin{l.}\label{dual}

If a totally real elliptic Lefschetz fibration admitting a real section is algebraic then the totally real elliptic Lefschetz fibration whose necklace diagram is dual to the necklace diagram of the former is also algebraic.\end{l.}

\begin{figure}[ht]
\begin{center}
      \includegraphics[scale=0.25,trim=0 0 70 0]{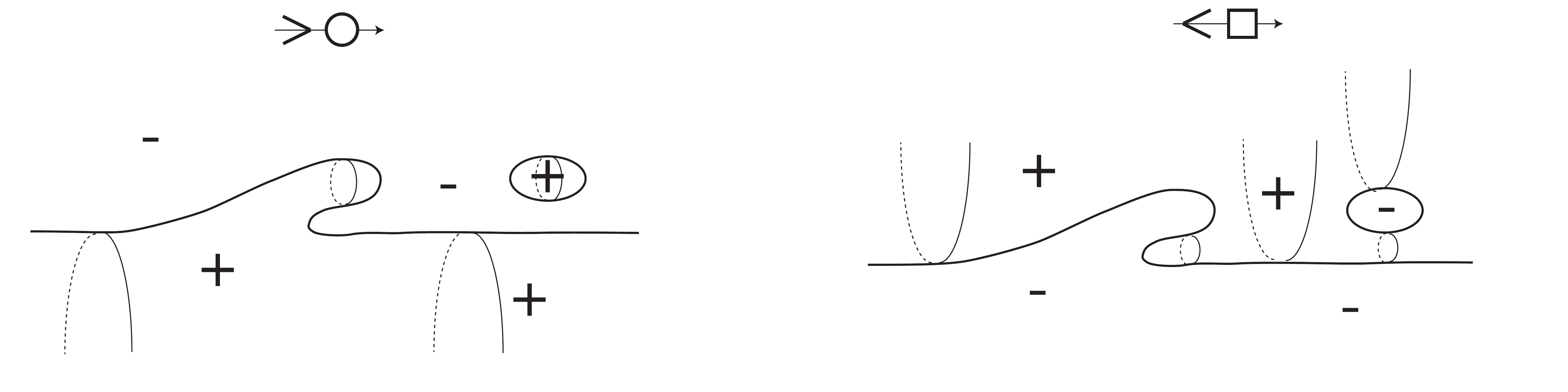}
\caption{\small{For each trigonal curve on $H(2n)$, there are two real structures of $E(n)$.}}
\label{Apisaret}
\end{center}
\end{figure}

\noindent {\it Proof: } 
The crucial observation is that although the real parts of fibrations associated with dual necklace diagrams  are  topologically different,  trigonal curves appearing as the branching set of coverings $E(n)\to H(2n)$ are the same. Duality of necklace diagrams corresponds, indeed, to the two different liftings of the real structure of $H(2n)$ to $E(n)$, see Figure~\ref{Apisaret}. 
 \hfill  $\Box$ \\

\begin{t.}\label{ceb} All totally real $E(1)$ admitting a real section are algebraic except those fibrations whose necklace diagram is one of the diagrams  listed in Figure~\ref{Apcebdegil}.
\end{t.}

\begin{figure}[ht]
\vspace{0.8cm}
\begin{center}
        \includegraphics[scale=0.19,trim=0 0 00 0]{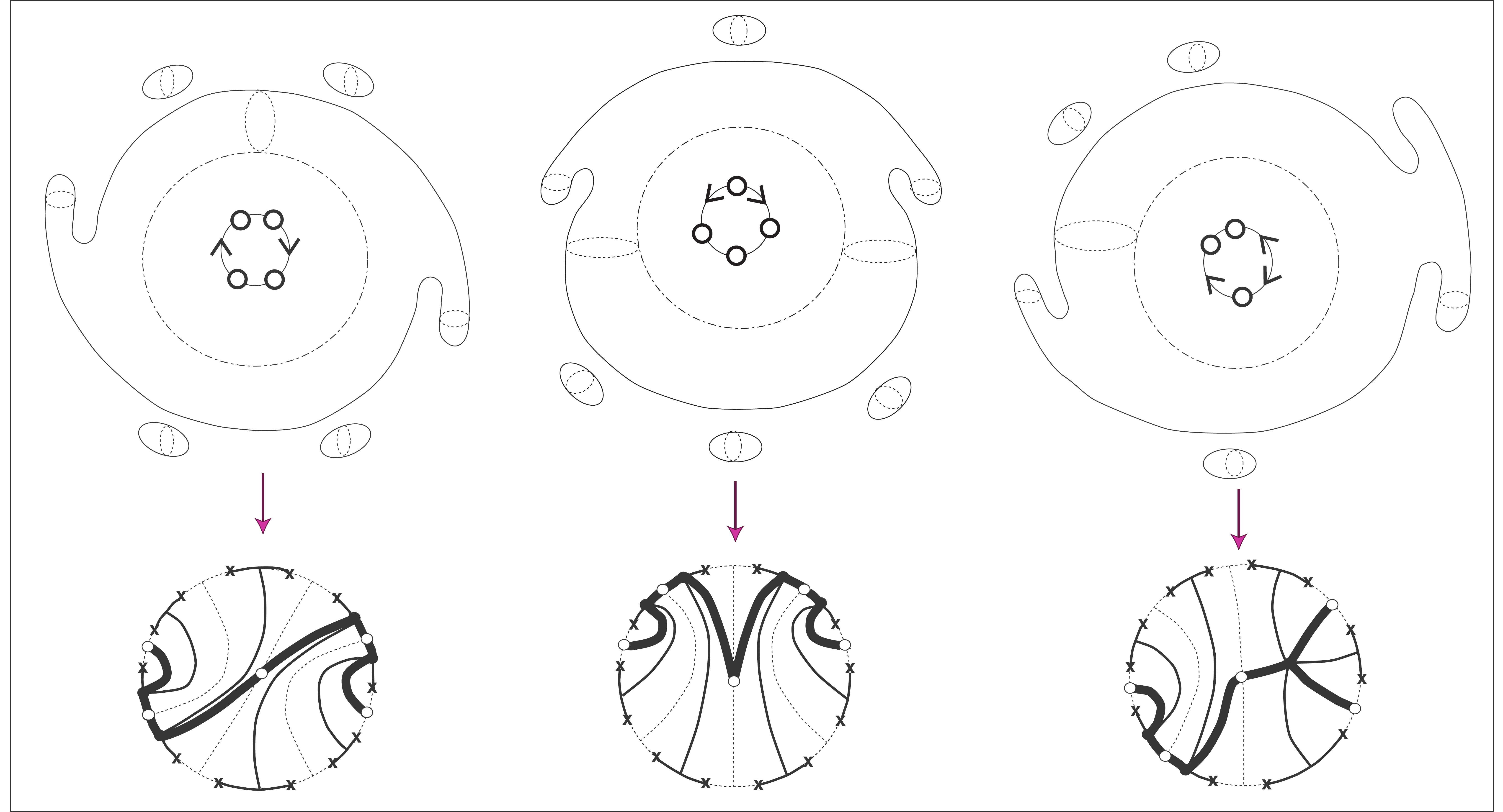}
 \caption{\small{Real dessin d'enfants.}}
\label{APalg1}
\end{center}
\end{figure}

\begin{figure}[ht]
\begin{center}
       \includegraphics[scale=0.19,trim=0 0 00 0]{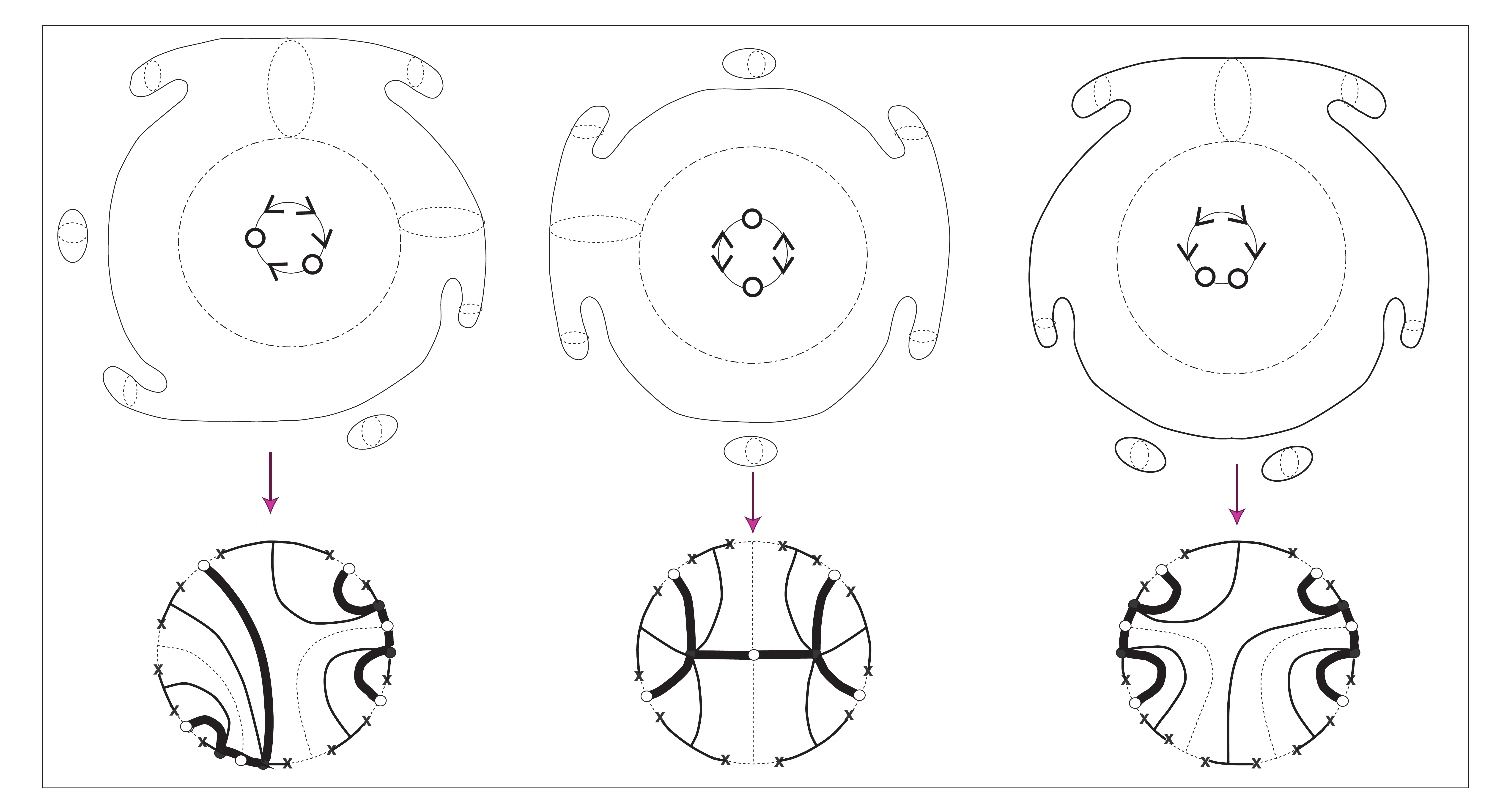}
\caption{\small{Real dessin d'enfants.}}
\label{APalg2}
\end{center}
\end{figure}

\begin{figure}[ht]
\begin{center}
        \includegraphics[scale=0.19,trim=0 0 00 0]{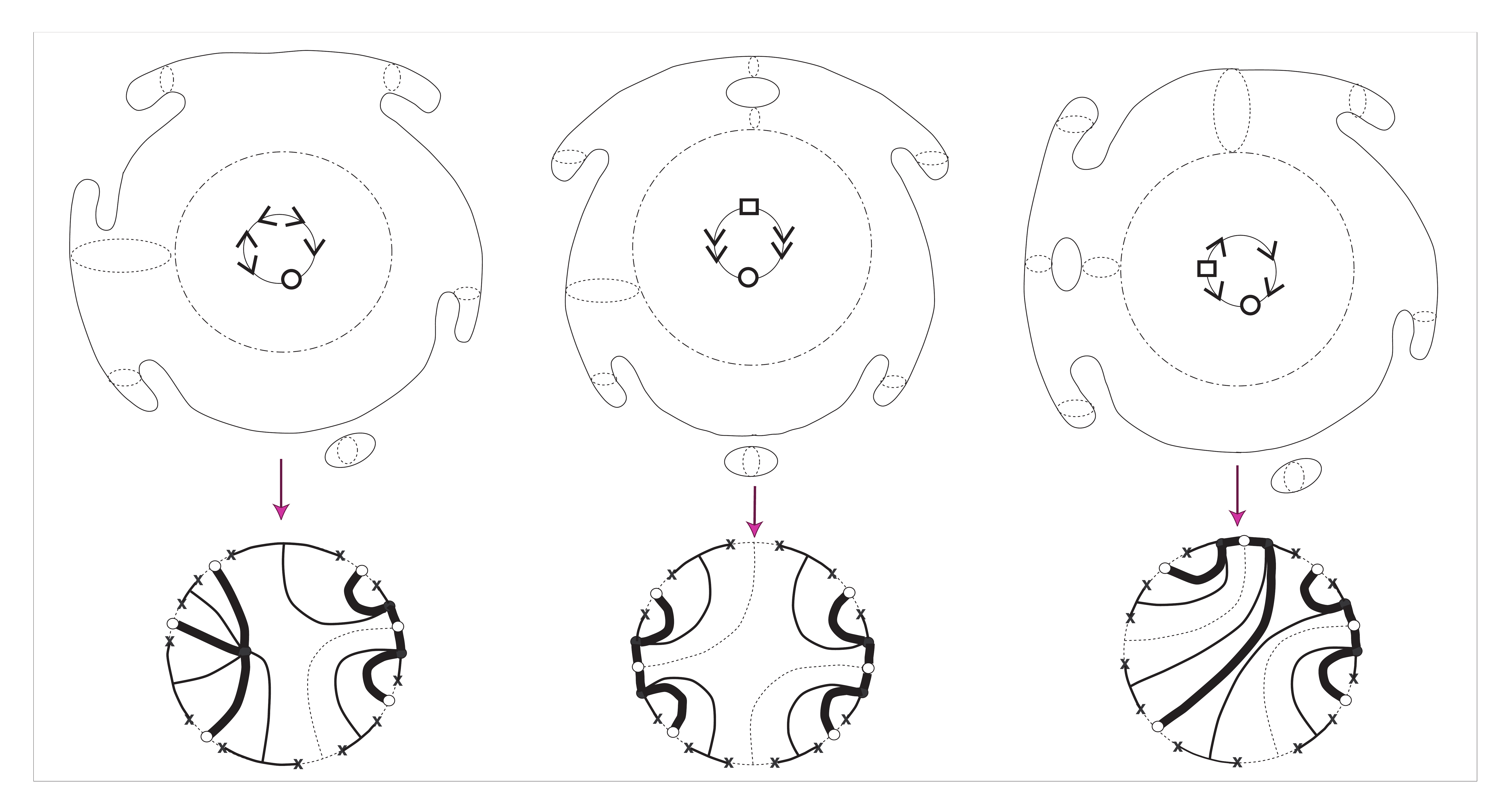}
        \caption{\small{Real dessin d'enfants.}}
\label{APalg3}
\end{center}
\end{figure}

\begin{figure}[ht]
\begin{center}
        \includegraphics[scale=0.19,trim=0 0 00 0]{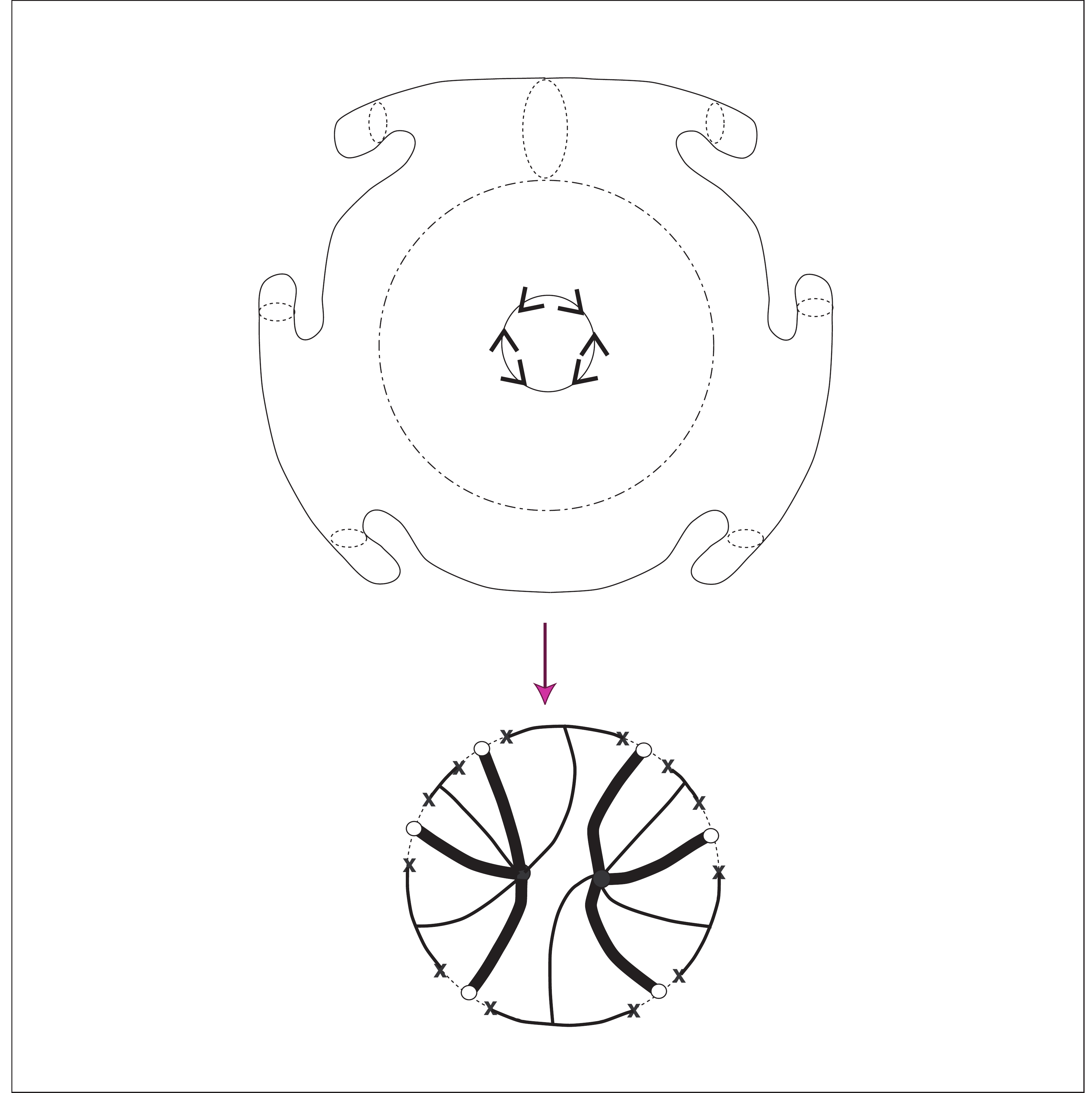} 
\caption{\small{Real dessin d'enfants.}}
\label{APalg4}
\end{center}
\end{figure}

\noindent {\it Proof: }  
By Theorem~\ref{torev} it is enough to construct real dessins d'enfants corresponding to necklace diagrams which are
not prohibited by Proposition~\ref{yasak}.
Following  Lemma~\ref{dual}, we only need to consider necklace diagrams with $|\ov|\geq |\kar|$.  
Figures~\ref{APalg1}-\ref{APalg4} show the required  real dessins d'enfants. \hfill  $\Box$ \\

\textbf{Real dessins d'enfants of real algebraic $E(1)$ with real sections.} (Around necklace diagrams, the real part is depicted. The dotted inner circle stands for  the lift of the exceptional section. Because of the symmetry, we only draw a half of the graph.)
        


\section{Necklace calculus and further applications}

In this section, we consider certain operations on the set of necklace diagrams.  These operations allow us to construct new necklace diagrams from the given ones.

\subsection{Necklace sums}

\begin{figure}[ht]
\begin{center}
        \includegraphics[trim = 5 5 0 0, clip, width=6cm]{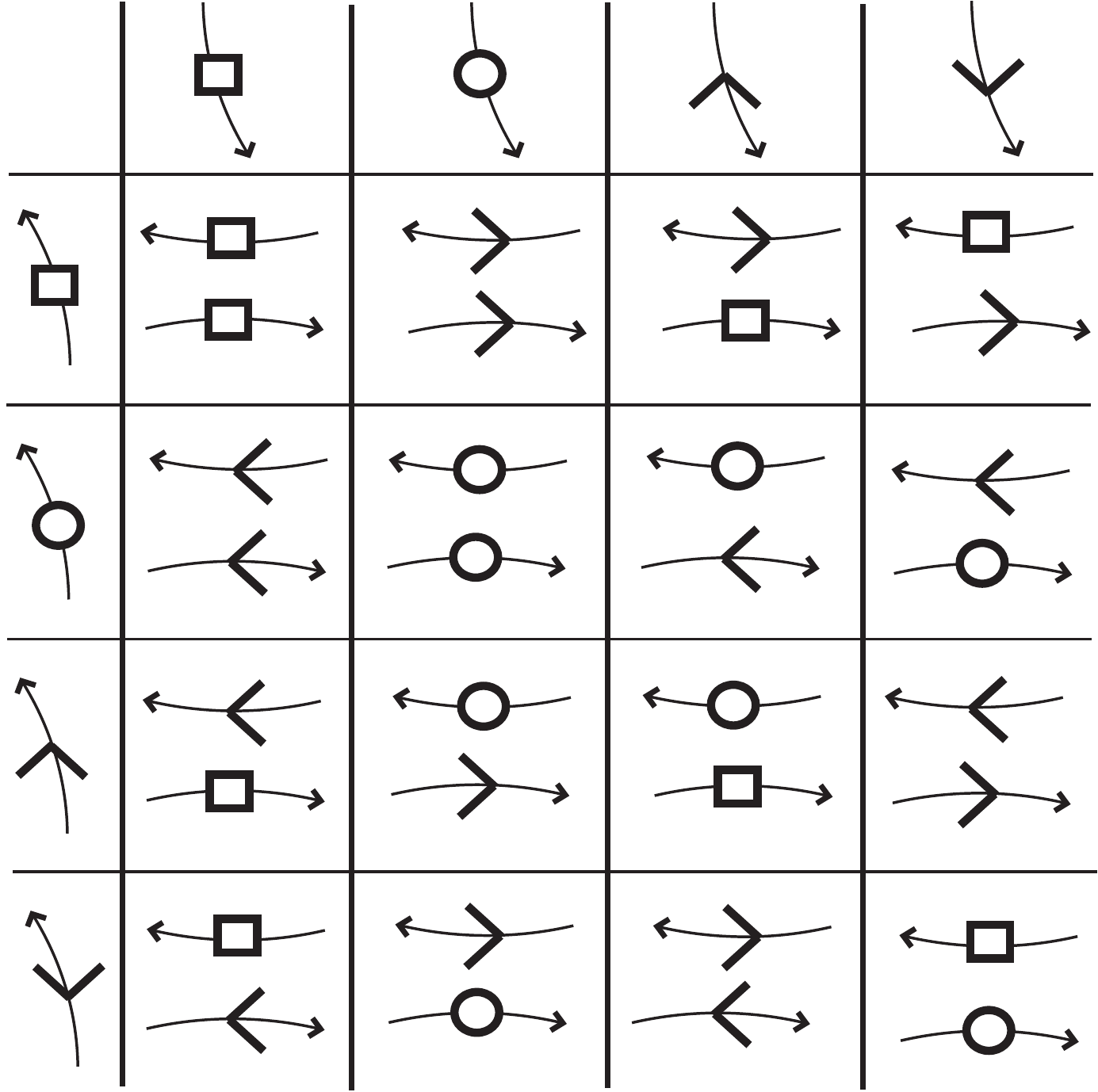}
             \caption{\small{Table of the harsh sum}}
             \label{tablo}
             \end{center}
 \end{figure}

A \emph{necklace sum} is basically the connected sum of the underlying oriented circle and it refers to the fiber sum of  the corresponding real Lefschetz fibrations. 
We consider two types of necklace sums which we call \emph{mild sum} and \emph{harsh sum}.
To perform a mild sum, we cut  each necklace diagram at a point on the chain then reglue the diagrams crosswise respecting the orientation.
The harsh sum, on the other hand, is obtained by cutting necklace diagrams at a stone and regluing them according to the table shown in Figure~\ref{tablo}. 
It follows from their definition that both the mild sum and the harsh sum do not change the monodromy of the diagram.
Evidently, the Euler characteristic is additive with respect to the mild sum;  however, it is not always additive with respect to the harsh sum. 

Let us note also that we can also consider necklace sum of non-oriented necklace diagram by fixing auxiliary orientations on diagrams.

\begin{figure}[ht]
   \begin{center}
     \includegraphics[trim = 5 5 0 0, clip, width=9cm]{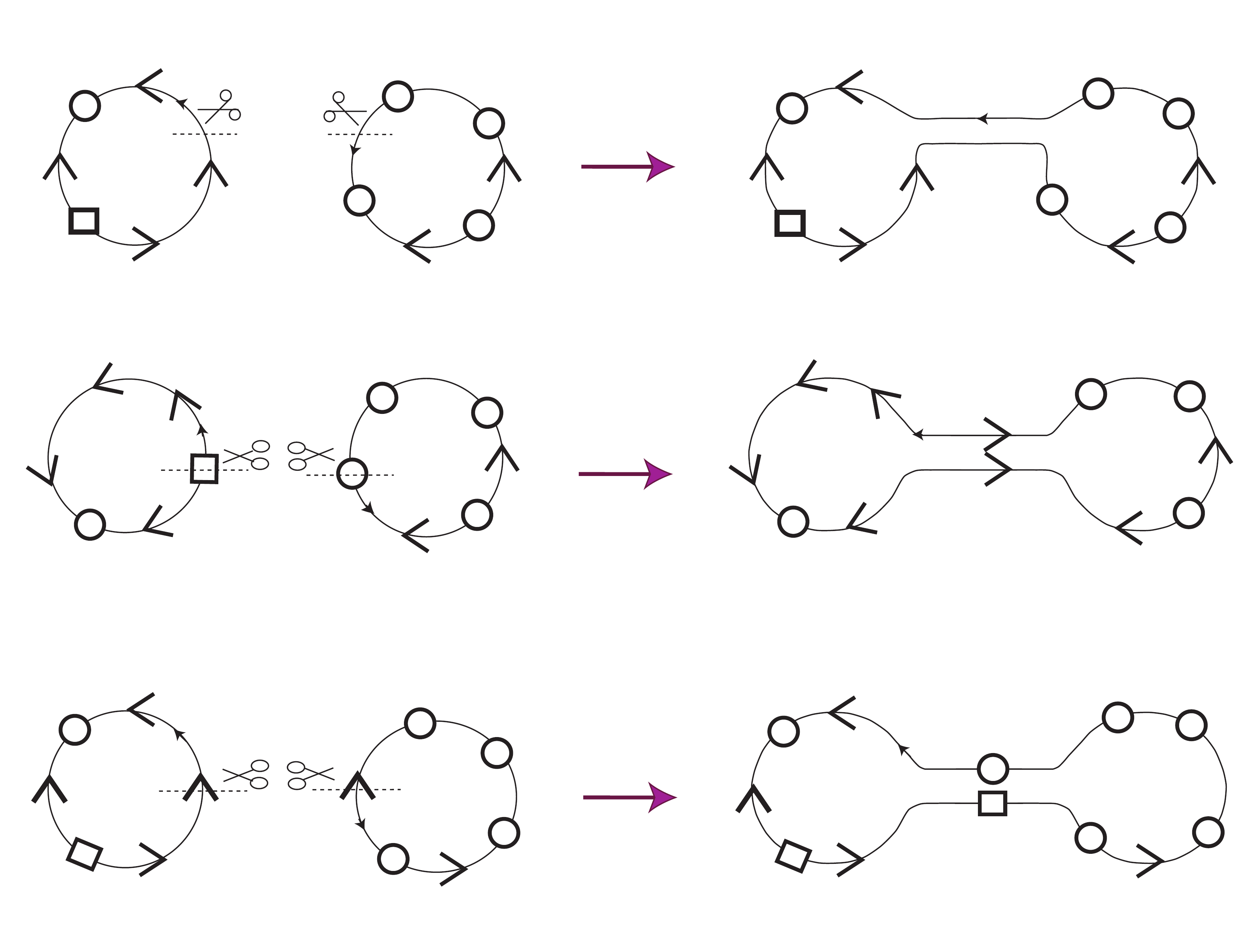}
       \caption{\small{Examples of the mild and the harsh sums.}}
       \label{ornek}
       \end{center}
 \end{figure}

Examples of mild and harsh sums are given in Figure~\ref{ornek}.

\begin{r.} There are two types of necklace chain segments (essential, non-essential) distinguished by the associated  graph fragments.
It is not hard to see that the mild sum preserves algebraicity if the points where the sum is taken are chosen on the same type of  chain segments and if the segments after the sum remain of the same type; or if they are chosen on different types of chain segments.  (In other words,  algebraicity is preserved  if we make a mild sum at two essential (respectively non-essential) intervals, in a way that after gluing we obtain again two essential (respectively non-essential) intervals; or if we make sum at an essential and a non-essential intervals.)  As for the harsh sum, we note that it   preserves algebraicity  if the number of neither $\ov$-type nor $\square$-type stone decreases after the sum. 
\end{r.}

\begin{p.}
For each $n$, maximal necklace diagrams exist and 
each totally real elliptic $E(n)$ represented by a maximal necklace diagram is algebraic.
\end{p.}

\noindent {\it Proof: } 
It is easy to see that the harsh sum of two maximal necklace diagrams where the sum is performed at arrow type stones of the opposite directions is maximal.
Moreover, by the remark above harsh sum performed at two arrow type stones of opposite directions preserves algebraicity. 
Note also that  as Theorem~\ref{ceb} asserted all maximal necklace diagrams of 6 stones are algebraic.

To finish the proof,  we show that  all maximal necklace diagrams of $n$ stones are obtained as harsh sums of maximal necklace diagrams of 6 stones.
 We will prove the claim by induction on $n$.
The first step is to check the claim for  $n=2$.  As we have the explicit list of maximal diagrams of 12 stones, we see immediately that any maximal necklace diagram of 12 stones can be obtained as the harsh sum of maximal necklace diagrams of 6 stones.  
Now, let us assume that for $n=k$ the claim is true. 
To prove the claim for  $n=k+1$,  note that the monodromies of $\ov$-type stones and $\kar$-type stones do not have any cancelation. 
The fact that there is no cancellation between $\ov$-type stones and $\kar$-type stones and that the monodromy of the necklace diagram is the identity
impose certain conditions on the possible arrangements of stones around an arrow type stone on a maximal necklace diagram. 
By checking the possibilities of the neighborhood of an arrow type stone, we see that required cancelations appear only in the cases where the arrangements come from the maximal necklace diagrams of 6 stones, so the claim follows from the inductive step.
\hfill  $\Box$ \\

\subsection{Flip-flops and metamorphoses}

Let ${\mathcal{N}}_k$ denote the set of (oriented) necklace diagrams with  $k$ stones and with the identity monodromy, and  let ${\mathcal{N}}_k^{(i, j)}$ be the subset consisting of diagrams with $(|\ov|, |\kar|)=(i, j)$. We define two kinds of operations on ${\mathcal{N}}_k$:   \emph{flip-flop} and  \emph{metamorphosis}. The  flip-flop preserves  $(|\ov|, |\kar|)$  and it coincides with canceling and creating handles on the real part  $E(n)_{\R}$. On the other hand, the metamorphosis decreases $|\ov|$ or $|\kar|$ by one and it appertains to a nodal deformation of $E(n)_{\R}$.  

\begin{figure}[ht]
   \begin{center}
        \includegraphics[scale=0.26,trim=0 0 0 0]{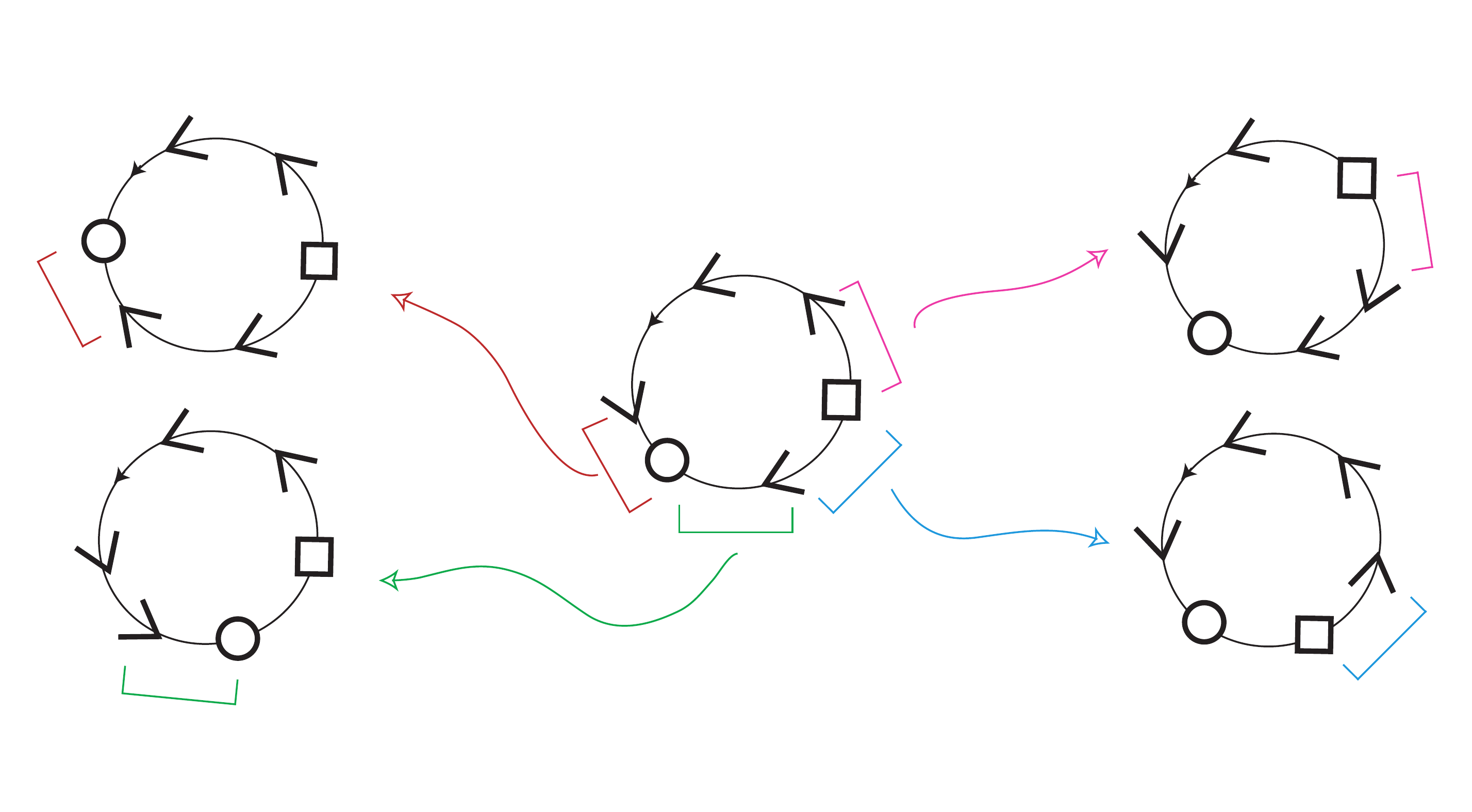}
           \caption{\small{Examples of flip-flops.}}
           \label{flipflopor}
\end{center}
\end{figure}

Flip-flop is the operation which swaps the segments shown below. Because the segments have the same monodromy, the total monodromy does not change. 
Examples of flip-flop are shown in Figure~\ref{flipflopor}.

\begin{displaymath}
\begin{array}{cccc}
\mbox{flip-flop}:& {{\mathcal{N}}_k}^{(i,j)}& \leftrightarrow & {{\mathcal{N}}_k}^{(i,j)} \\

     &\includegraphics[trim= 0 0 0 0, clip, width=1.3cm]{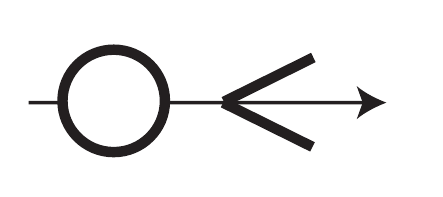}&  \includegraphics[trim = 0 0 0 0, clip, width=1cm]{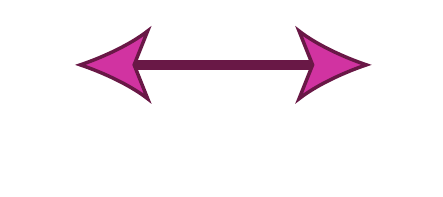}    & \includegraphics[trim = 0 0 0 0, clip, width=1.3cm]{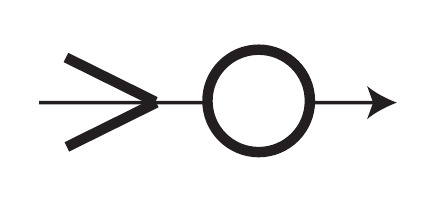}\\

     &\includegraphics[trim = 0 0 0 0, clip, width=1.3cm]{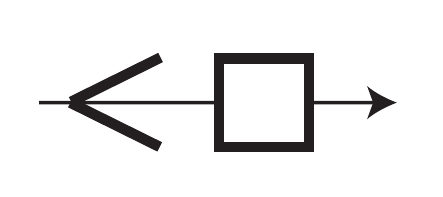}&    \includegraphics[trim = 0 0 0 0, clip, width=1cm]{cift.pdf}  & \includegraphics[trim = 0 0 0 0, clip, width=1.3cm]{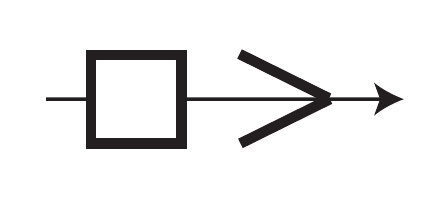}
\end{array}
\end{displaymath}

\begin{figure}[ht]
   \begin{center}
       \includegraphics[scale=0.24,trim=0 0 0 0]{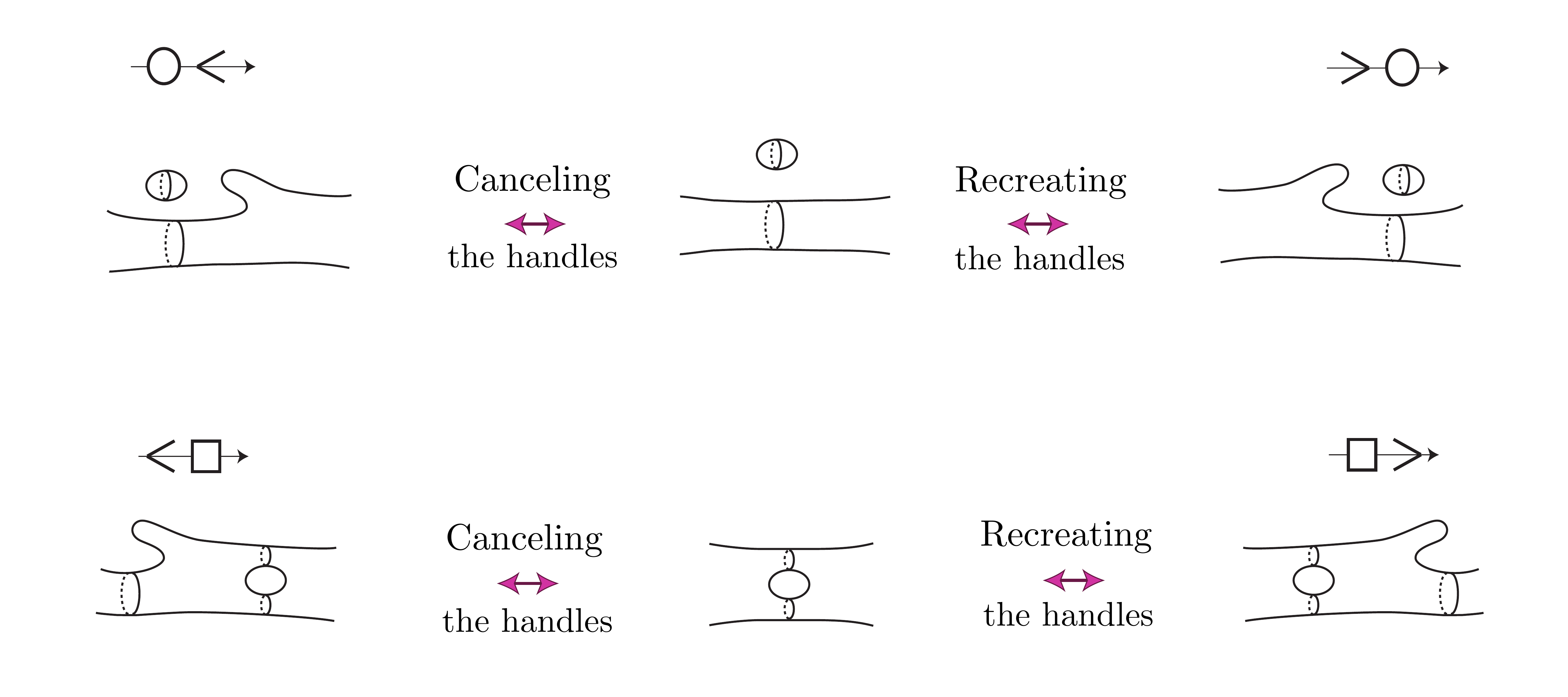}
         \caption{\small{The effect of  flip-flop on  the real part.}}
         \label{flipflopreel}
\end{center}
\end{figure}

As $|\ov|$ and $|\kar|$ remain unchanged after a flip-flop, the Euler characteristic and the total Betti number of  the corresponding $E(n)_{\R}$ do not change. Thus,  topological type of $E(n)_{\R}$ is not affected by the flip-flop. In Figure~\ref{flipflopreel}, we interpret the effect of a flip-flop on $E(n)_{\R}$.

\begin{figure}[ht]
   \begin{center}
         \includegraphics[scale=0.26,trim=0 0 70 0]{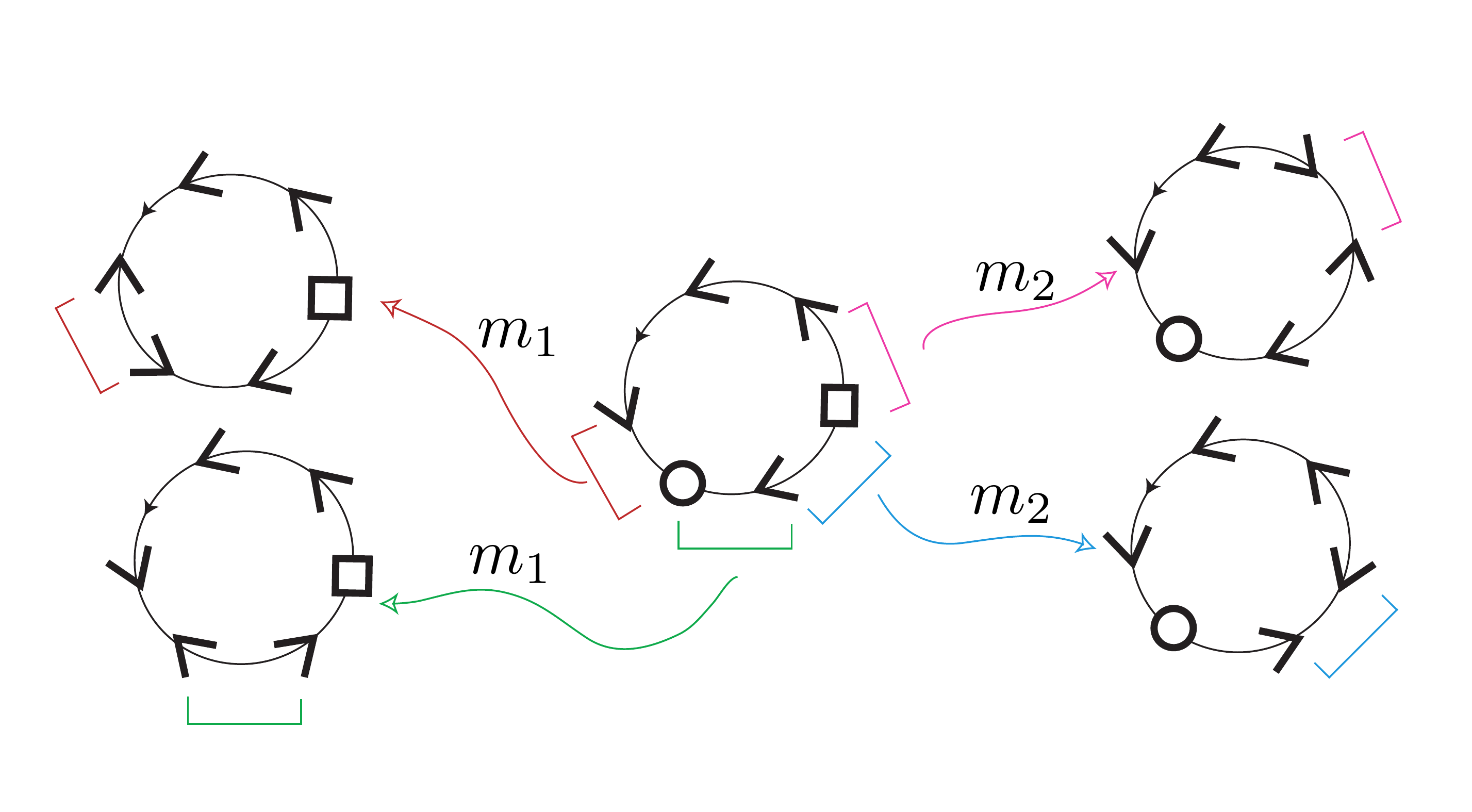}
          \caption{\small{Examples of metamorphoses.}}
\label{metaor}
\end{center}
\end{figure}

We consider two types of metamorphoses,  $m_1, m_2$, of necklace diagrams.
They modify the stones of the above segments  as below. 
Examples of metamorphoses are depicted in Figure~\ref{metaor}.

\begin{displaymath}
\begin{array}{cccc}
 m_1:& {{\mathcal{N}}_k}^{(i,j)}& \rightarrow & {{\mathcal{N}}_k}^{(i-1,j)} \\

     & \includegraphics[trim = 0 0 0 0, clip, width=1.3cm]{ox.pdf}& \includegraphics[trim = 0 0 0 0, clip, width=1cm]{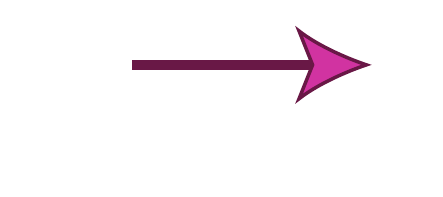}    &  \includegraphics[trim = 0 0 0 0, clip, width=1.3cm]{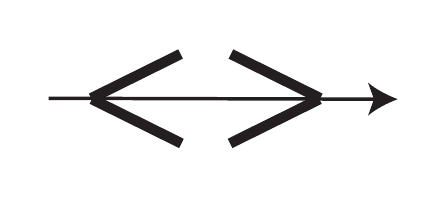}\\

     & \includegraphics[trim = 0 0 0 0, clip, width=1.3cm]{xo.pdf}&  \includegraphics[trim = 0 0 0 0, clip, width=1cm]{tek.pdf}   & \includegraphics[trim = 0 0 0 0, clip, width=1.3cm]{disa.pdf}

     \end{array}
\end{displaymath}

\begin{displaymath}
\begin{array}{cccc}
 m_2:& {{\mathcal{N}}_k}^{(i,j)}& \rightarrow & {{\mathcal{N}}_k}^{(i,j-1)} \\

     &\includegraphics [width=1.3cm]{ox2.pdf}& \includegraphics[trim = 0 0 0 0, clip, width=1cm]{tek.pdf}    & \includegraphics [width=1.3 cm]{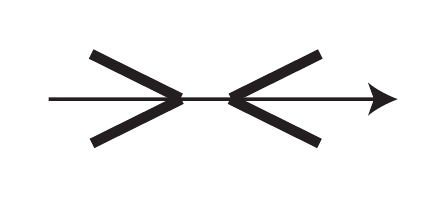}\\

     &\includegraphics [width=1.3cm]{xo2.pdf}&  \includegraphics[trim = 0 0 0 0, clip, width=1cm]{tek.pdf}   & \includegraphics [width=1.3 cm]{ice.pdf}
    \end{array}
\end{displaymath}

\begin{figure}[ht]
   \begin{center}
          \includegraphics[scale=0.19,trim=0 0 -40 0]{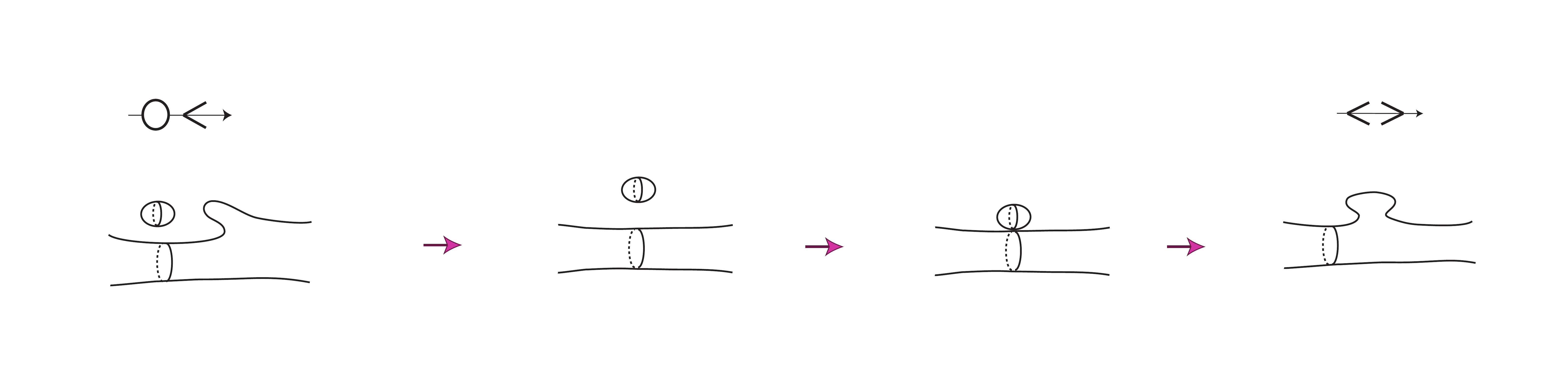}
                   \includegraphics[scale=0.19,trim=0 0 70 0]{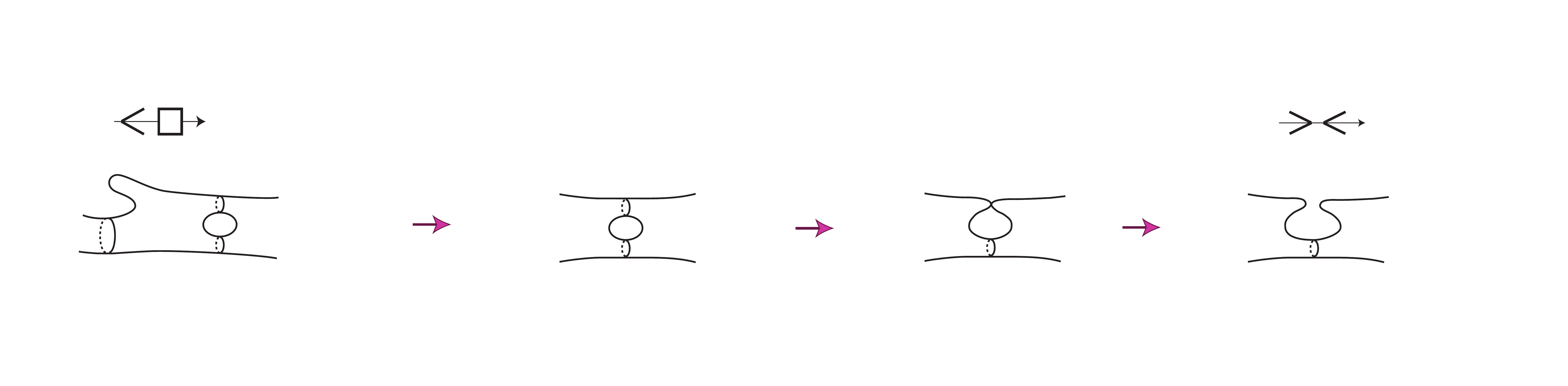}
                  \caption{\small{The effect of metamorphoses, $m_1, m_2$, to the real part.}}
         \label{metareel}
\end{center}
\end{figure}

Since $|\ov|$ or $|\kar|$ are modified by metamorphoses, the topological type of the corresponding $E(n)_{\R}$ changes.  Recall that for fibrations which admit a real section each $\ov$-type stone corresponds to a spherical component while each $\kar$-type stone drives a handle. Indeed, a sphere component  or a genus disappears after a metamorphosis (or appears after a inverse metamorphosis).  In Figure~\ref{metareel},  we  depict the effects of metamorphoses on $E(n)_{\R}$.

\subsection{Applications}

Below, in Figure~\ref{e1top} and Figure~\ref{e2top}, we present two graphs (for $E(1)$ and $E(2)$, respectively) whose vertices correspond to necklace diagrams with fixed $(|\ov|, |\kar|)$, edges to necklace metamorphoses  $m_{1}, m_{2}$.  As we mention before the real part of totally real  $E(n)$,  admitting a real section,  consists of  spherical components (the number of which is $|\ov|$) and a higher genus component which is an orientable surface of genus  $|\kar|+1$  if $n$ is even; a non-orientable surface with $2|\kar|+1$ cross-caps, otherwise.  Each pair $(|\ov|, |\kar|)$ and the parity of $n$, thus, defines  the topological type of $E(n)_{\R}$.

\begin{figure}[ht]
\begin{center}
   \includegraphics[trim = 5 5 0 0, clip, width=6cm]{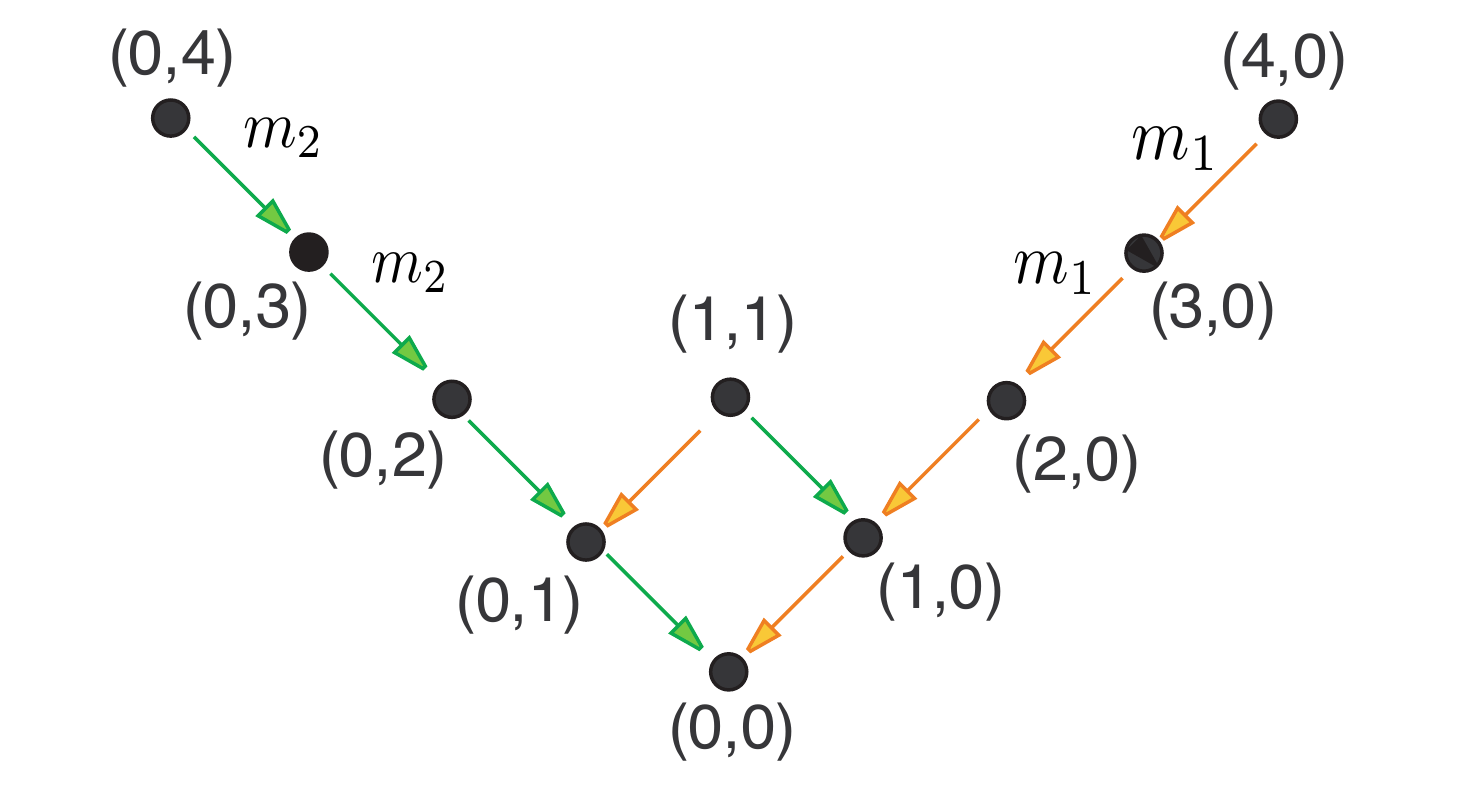}
      \caption{\small{ Metamorphosis graph of  $E(1)_{\R}$.}}
       \label{e1top}
\end{center}
\end{figure}

\begin{figure}[ht]
\begin{center}
   \includegraphics[trim = 5 5 0 0, clip, width=13cm]{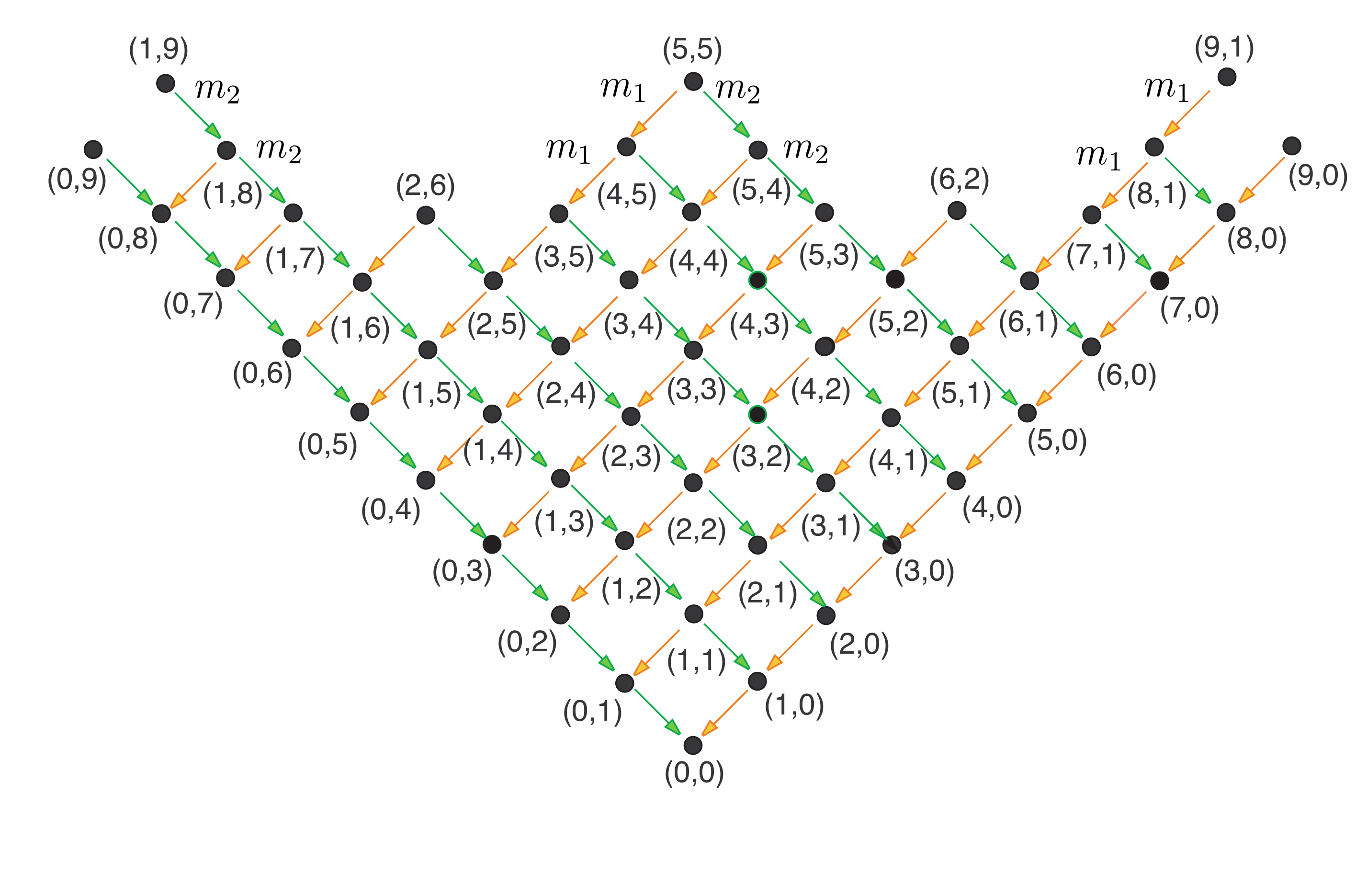}
       \caption{\small{ Metamorphosis graph of $E(2)_{\R}$. }}
       \label{e2top}
       \end{center}
 \end{figure}

Examining the list of necklace diagrams of 6 stones we obtain the following:

 \begin{p.}
All 6-stone necklace diagrams  listed in Figure~\ref{e1kolyeler}  can be obtained from the maximal ones by a sequences of metamorphoses, inverse metamorphoses or flip-flops.
Moreover, the list of necklace diagrams which are obtained from maximal diagrams only by a sequence of metamorphoses and eventually an inverse metamorphosis, coincide with the list of diagrams of algebraic fibrations. 
  \hfill  $\Box$ \\
\end{p.}

\begin{p.} There exist 12-stone necklace diagrams (with the identity monodromy) that can not be obtained from
the maximal necklace diagrams by necklace operations.  
\end{p.}

\begin{figure}[h]
\begin{center}
      \includegraphics[scale=0.26,trim=0 0 70 0]{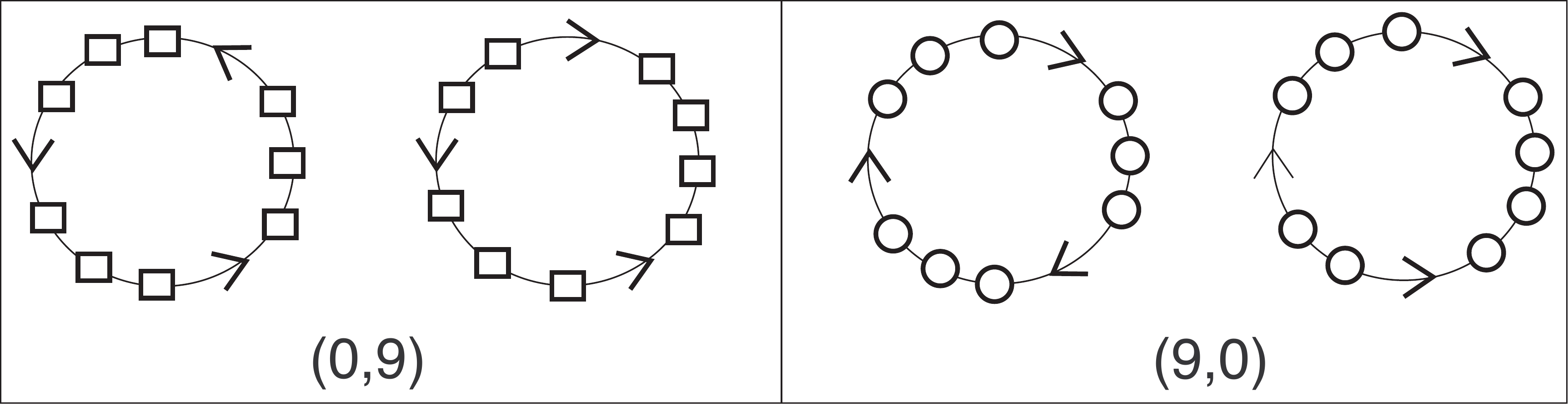}
       \caption{\small{Necklace diagrams which can not be obtained by necklace operations.}}
       \label{90ex}
       \end{center}
 \end{figure}

\noindent {\it Proof: } 
By duality it is enough to consider the case of $|\ov|\geq |\kar|$.
Examples, shown in Figure~\ref{90ex},  are found by investigating the list of 12-stone necklace diagrams  basically with  $(|\ov|,|\kar|)=(9,1), (9,0), (8,1), (8,0)$ .
 Let us also note that the  fibrations associated  with the diagrams  shown in Figure~\ref{90ex}  are algebraic.
\hfill  $\Box$ \\

  As a corollary of the above proposition we claim that all totally  real algebraic  $E(1)$ admitting a real section can be obtained from the maximal ones by a sequences of nodal deformations,  while there are totally real algebraic  $E(2)$ which cannot be obtained in this way.

\begin{p.}
There exist 12-stone necklace diagrams (with the identity monodromy)  which 
 are not a necklace sum of two 6-stone necklace diagrams listed in \ref{e1kolyeler}.
\end{p.}

\begin{figure}[ht]
\begin{center}
     \includegraphics[scale=0.22,trim=0 0 10 0]{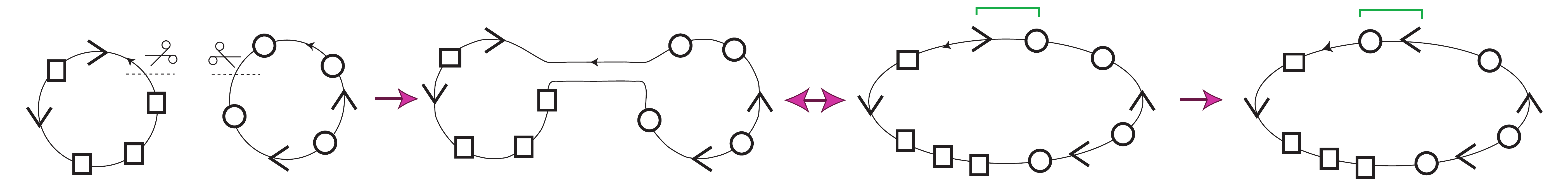}
       \caption{\small{An example of construction of a non-decomposable necklace diagram.}}
       \label{4notconnect}
       \end{center}
 \end{figure}

\noindent {\it Proof: } In Figure~\ref{4notconnect}, we construct a non-decompasable  example applying a mild sum followed by a flip-flop.
Let us also note that such examples can be produced the same way for any $n>1$.
By analyzing possible divisions of the pair $(|\ov|, |\kar|)$, we see that the necklace diagram shown in
Figure~\ref{4notconnect} cannot be divided into two 6-stone necklace
diagrams with the identity.
 \hfill  $\Box$ \\

\begin{c.}
There exist  totally real $E(2)$ which cannot be written as a fiber sum of two totally real $E(1)$.
 \hfill  $\Box$ \\
\end{c.}


\section{Totally real elliptic Lefschetz fibrations without a real section and refined necklace diagrams}
In this section, we explore  the case of totally  real elliptic Lefschetz fibrations $\pi:X\to S^2$ which do not admit a real section and introduce \emph{refined necklace diagrams}  associated with  them.

\begin{figure}[ht]
   \begin{center}
        \includegraphics[scale=0.3,trim=0 0 70 0]{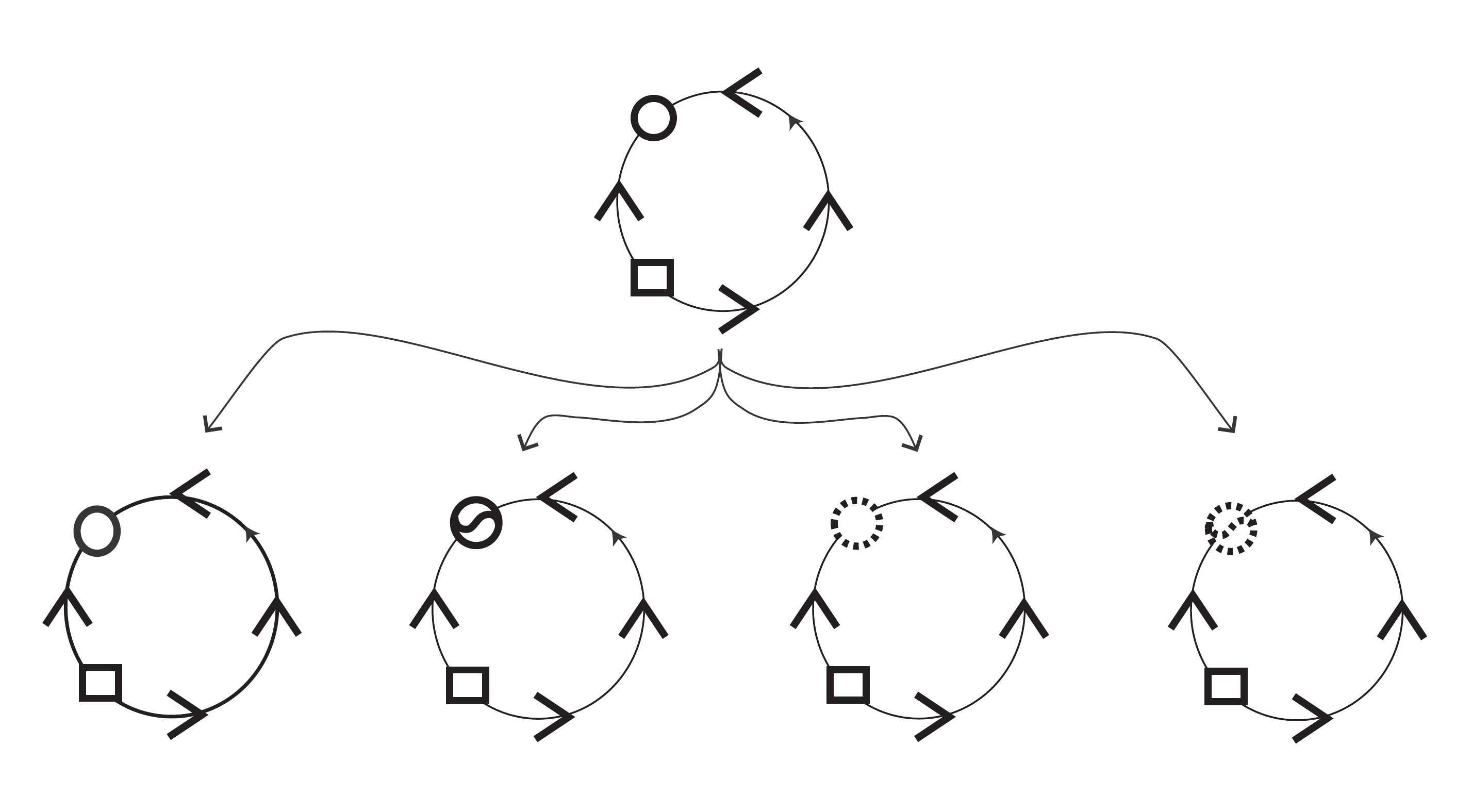}
\caption{\small{Refinement of a necklace diagram.}}
\label{incikolye}
       \end{center}
 \end{figure}

A refinement of a necklace diagram is obtained by replacing each $\ov$-type stone with one of the following refined stones, $\ovbir, \oviki, \ovuc, \ovdort$. 
If the refined necklace diagram is identical to the underlying necklace diagram then the corresponding real Lefschetz fibration admits a real section.
Examples of refinements of a necklace diagram are shown in Figure~\ref{incikolye}.

From Remark~\ref{vanissaykil}, it is clear that  if a real structure on a fiber of $\pi$ has no real component, then the nearby critical values can only be of type ``$\circ$".  In other words,  existence or lack of a real section influences only  $\ov$-type necklace stones.

\begin{figure}[ht]
   \begin{center}
     \includegraphics[scale=0.3,trim=0 0 70 0]{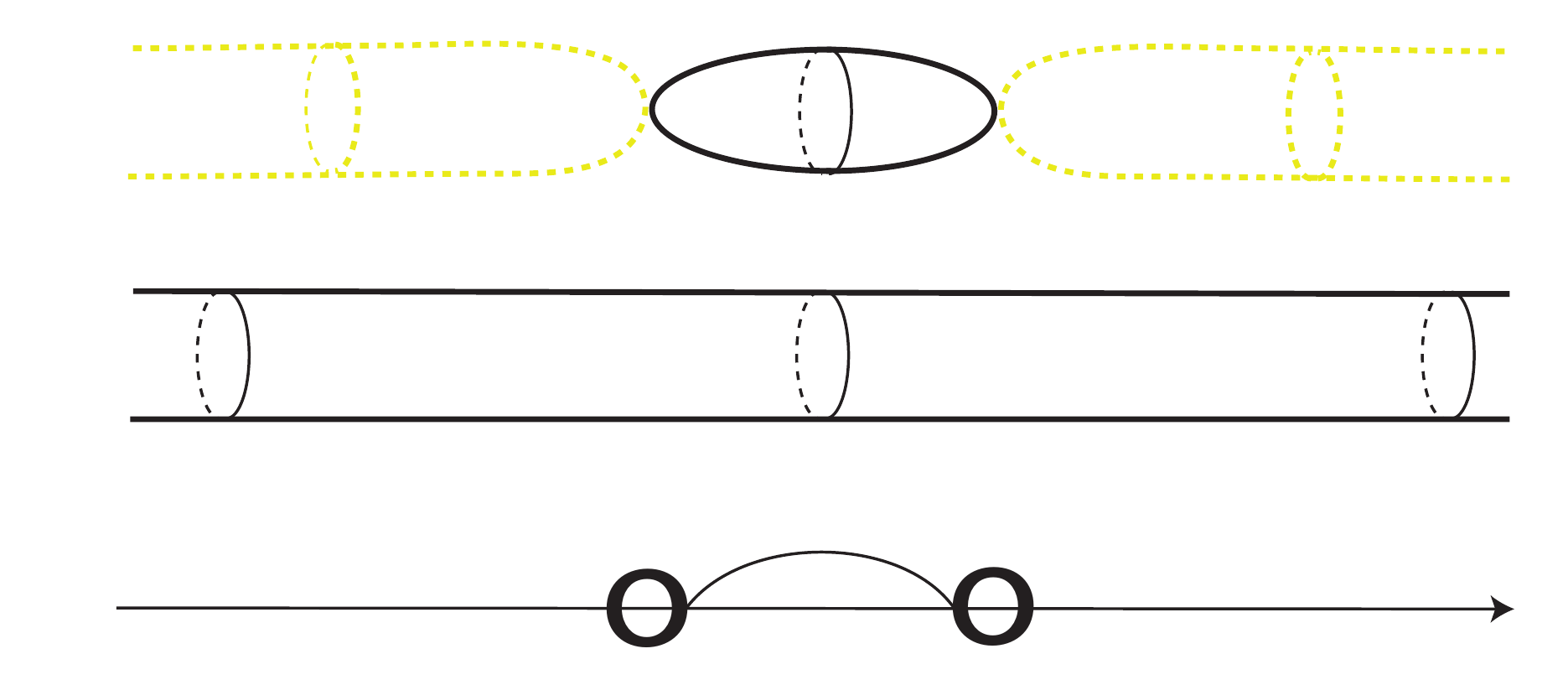}\hspace{1cm}
     \includegraphics[scale=0.3,trim=0 0 70 0]{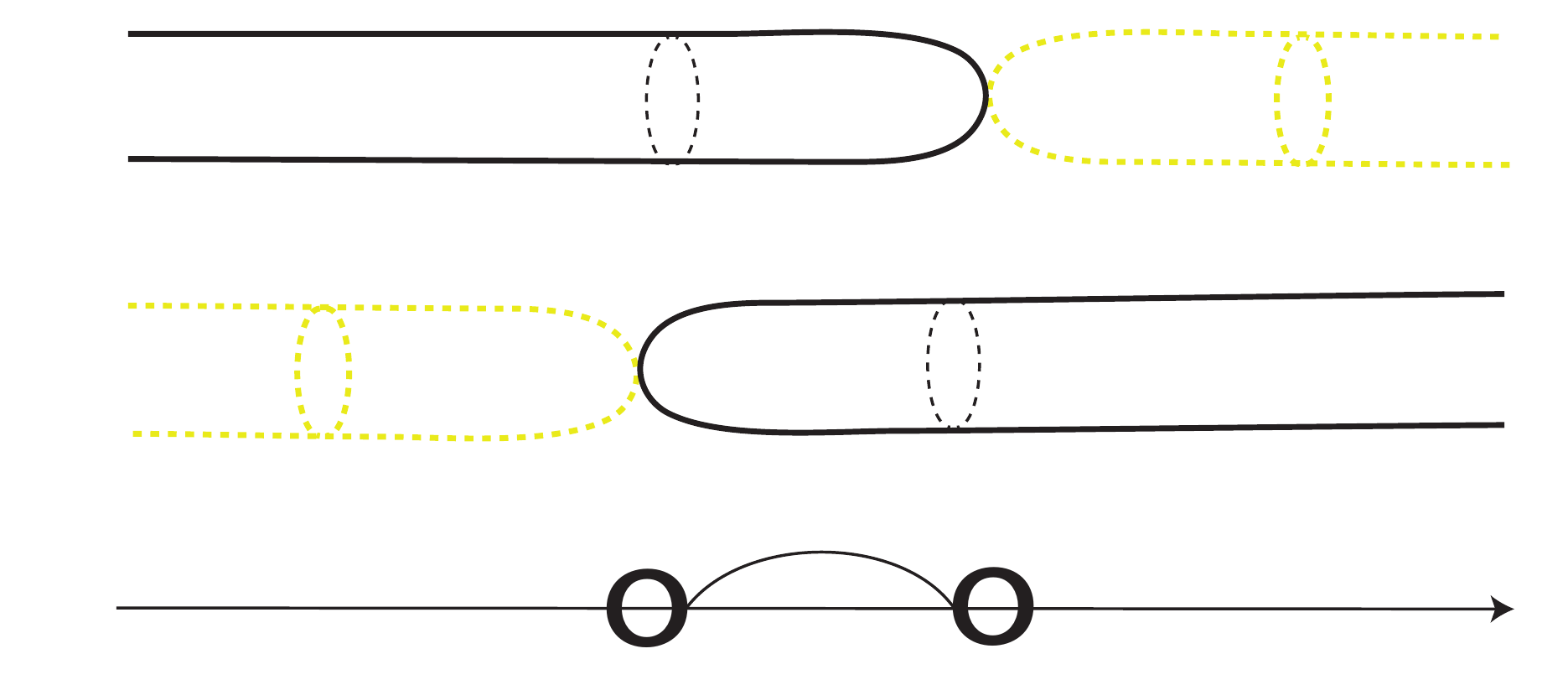}
      \caption{\small{Real part and associated refined stones.}}
      \label{sekt}
      \end{center}
 \end{figure}

Both of the refined stones of type $\ovbir, \oviki$ correspond to the case where the real structure on the real fibers over the interval between the two critical values has 2 real components. 
Already the real part $X_{\R}$ of $X$  distinguishes  the cases of $\ovbir$ and $ \oviki$, see Figure~\ref{sekt}.  As notation suggested $\ovbir$ has to do with the case where there is a real section, and hence, only  $\ovbir$-type refined stones refer to a spherical component of $X_{\R}$. 

\begin{figure}[ht]
   \begin{center}
         \includegraphics[scale=0.3,trim=0 0 70 0]{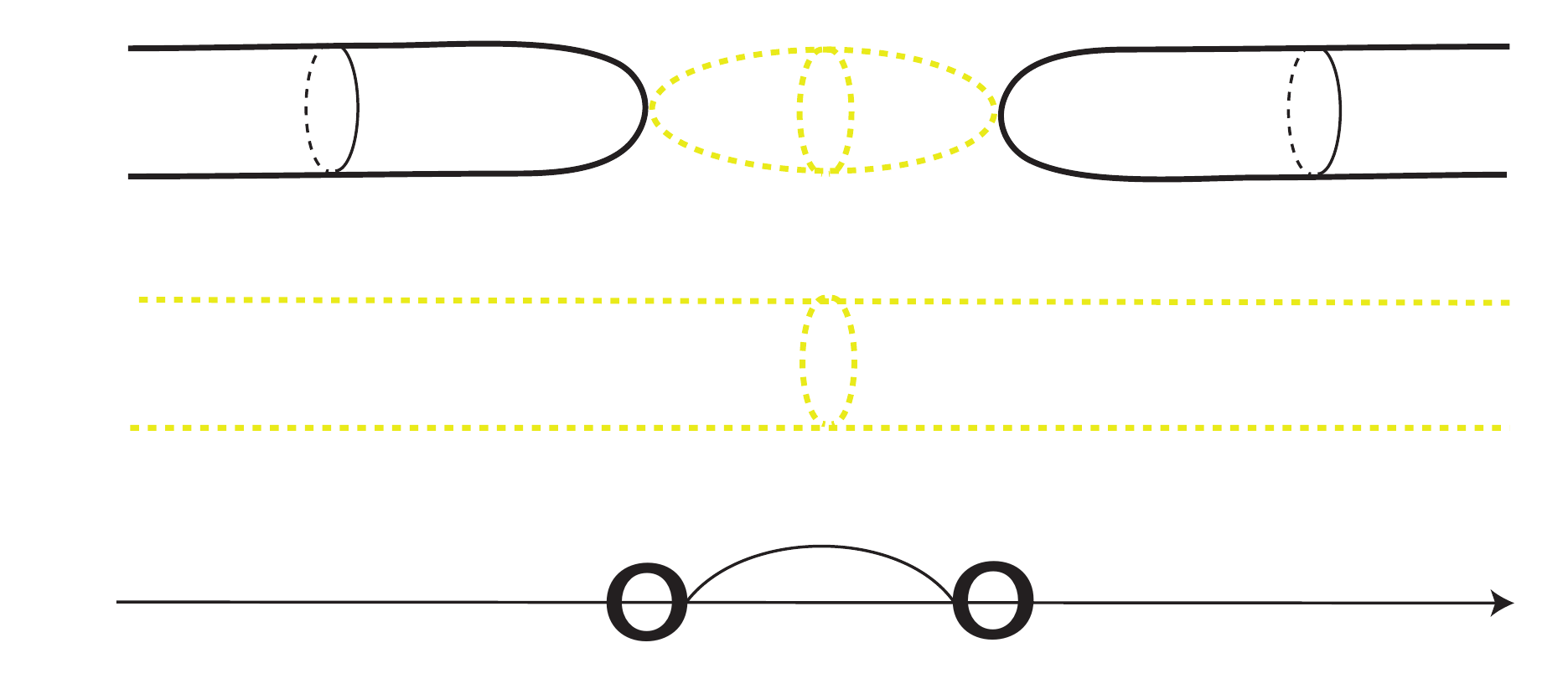}\hspace{1cm}
      \includegraphics[scale=0.3,trim=0 0 70 0]{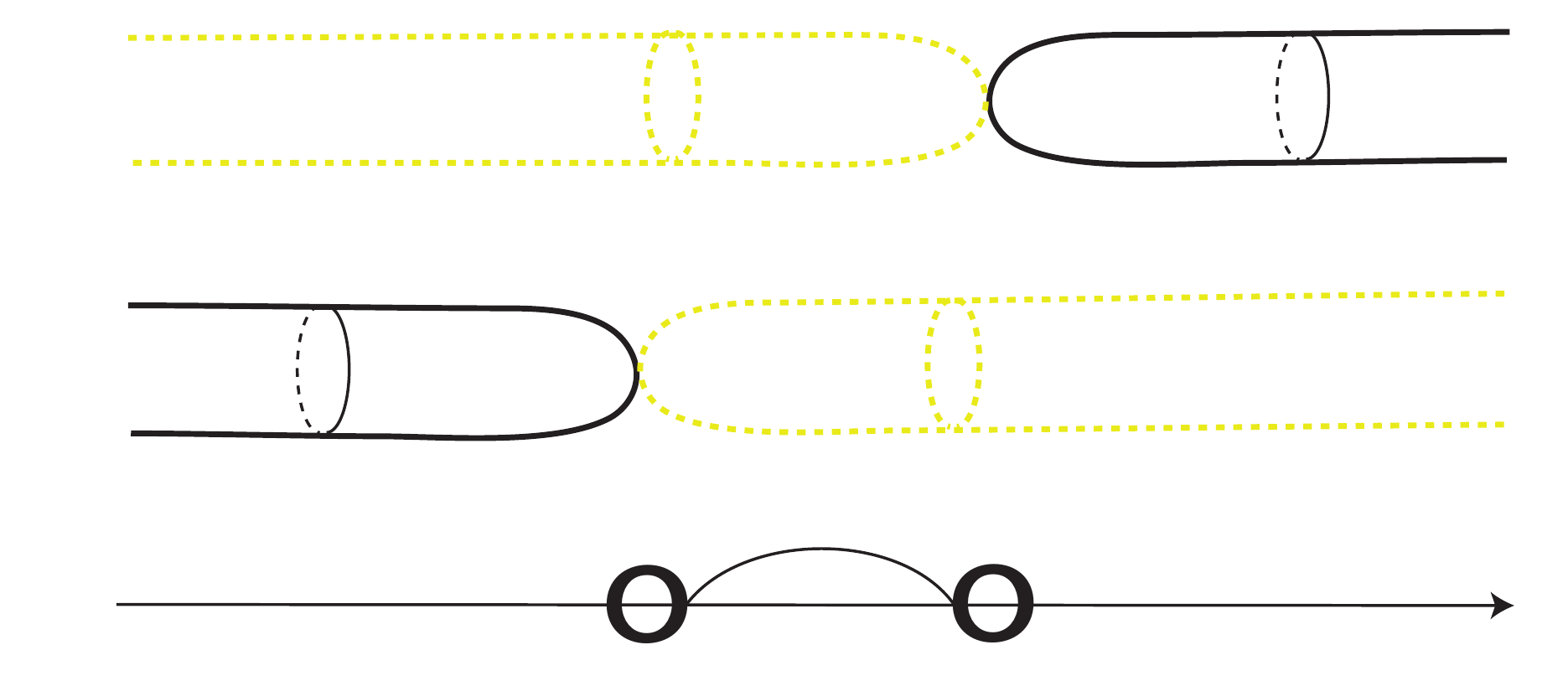}
        \caption{\small{Real part and associated refined stones.}}
     \label{ssekt}
       \end{center}
 \end{figure}

Recall that if the condition that the fibration admits a real section is discarded,  then the fibers of $\pi_{\R}$ may also be empty (which happens when the real structure on the real fibers of $\pi$ has no real component).  We introduce the refined stones  $\ovuc$ and $\ovdort$ which  correspond to the case where the real structure on the real fibers of $\pi$ has no real component. 
As depicted in Figure~\ref{ssekt}, the real part of $X_{\R}$  does not distinguish the two situations associated  with  $\ovuc, \ovdort$. 
The difference between $\ovuc$ and $\ovdort$ (as well as between $\ovbir$ and $\oviki$) can indeed be conceived by comparing the equivariant isotopy classes of  
the two vanishing cycles corresponding to the two critical values of the necklace stone.  In the case of $\ovuc$ (respectively $\ovbir$) the equivariant isotopy classes of vanishing cycles are the same, while in the case of $\ovdort$ (respectively, $\oviki$)
the two vanishing cycles  are of different  equivariant classes.  (A more detailed discussion can be found in  \cite[Section~8]{ner2}.)

As mentioned in Remark~\ref{mono}, there is no difference between real
structures with 2 real components and  real structures with no
component  on the homological level.  As a consequence, the calculation
of the monodromy does not affected by the refinement.  Thus, we have the following theorem.

\begin{t.}\label{birebirincikolye}
There is a one-to-one correspondence between the set of oriented refined necklace diagrams with $6n$ stones whose monodromy is the
identity and the set of isomorphism classes of directed totally real $E(n)$, $n\in\N$. 
\end{t.}

\noindent {\it Proof: }  The proof is analogous to the proof of Theorem~\ref{birebirkolye}.  
It is obvious from its construction that the refinements of necklace diagram is exactly the decoration of the real Lefschetz chains, see Figrure~10 of \cite{ner2}. Thus, we relate refined necklace diagrams with the  \emph{decorated real Lefschetz fibrations}, complete invariants of totally real elliptic Lefschetz fibrations, presented in  \cite[Section~8]{ner2}.
The result, thus, follows from Theorem~8.1 and  Proposition~8.2 of \cite{ner2}.  \hfill  $\Box$ \\

\begin{c.}\label{birebirsimliincikolyeler}
There is a one-to-one correspondence between the set of symmetry classes non-oriented refined necklace diagrams with $6n$ stones whose monodromy is the
identity and the set of isomorphism classes of  totally real $E(n)$, $n\in\N$. 
 \hfill  $\Box$ \\
\end{c.}

The number of possible refinements of necklace diagrams with fixed $(|\ov|, |\kar|)$ listed Figure~\ref{e1kolyeler} is given below. 
\begin{itemize}
\item $(|\ov|,|\kar|)=(1,1)$ there are 12 refined necklace diagrams,
\item $(|\ov|,|\kar|)=(1,0)$ there are 8 refined necklace diagrams,
\item $(|\ov|,|\kar|)=(2,0)$ there are 46 refined necklace diagrams,
\item $(|\ov|,|\kar|)=(3,0)$ there are 84 refined necklace diagrams,
\item $(|\ov|,|\kar|)=(4,0)$ there are 251 refined necklace diagrams.
\end{itemize}


\end{document}